\documentclass[a4paper,english,oneside]{smfart}
\usepackage{amsmath,amssymb,amsthm}
\usepackage{xcolor,microtype}

\usepackage{enumitem}
\usepackage{ytableau}
\ytableausetup{boxsize=1.55em,aligntableaux=center}

\usepackage{mathrsfs}

\ExplSyntaxOn
\NewDocumentCommand{\definealphabet}{mmmm}
 {%
  \int_step_inline:nnn { `#3 } { `#4 }
   {
    \cs_new_protected:cpx { #1 \char_generate:nn { ##1 }{ 11 } }
     {
      \exp_not:N #2 { \char_generate:nn { ##1 } { 11 } }
     }
   }
 }
\ExplSyntaxOff
\definealphabet{bb}{\mathbb}{A}{Z}
\definealphabet{cal}{\mathcal}{A}{Z}
\definealphabet{frak}{\mathfrak}{A}{z}
\definealphabet{rm}{\mathrm}{A}{z}
\definealphabet{bf}{\mathbf}{A}{z}
\definealphabet{scr}{\mathscr}{A}{Z}
\definealphabet{tt}{\mathtt}{A}{Z}
\definealphabet{ul}{\underline}{A}{z}

\usepackage[pdfencoding=auto,psdextra]{hyperref}
\usepackage[nameinlink]{cleveref}
\Crefname{theorem}{Theorem}{Theorems}
\Crefname{conjecture}{Conjecture}{Conjectures}
\Crefname{lemma}{Lemma}{Lemmas}
\Crefname{definition}{Definition}{Definitions}
\Crefname{remark}{Remark}{Remarks}
\Crefname{question}{Question}{Questions}
\Crefname{proposition}{Proposition}{Propositions}
\Crefname{corollary}{Corollary}{Corollaries}
\Crefname{equation}{}{}
\Crefname{item}{}{}
\Crefname{example}{Example}{Examples}
\Crefname{proof}{Proof}{Proofs}
\crefname{introthm}{Theorem}{Theorems}
\Crefname{introthm}{Theorem}{Theorems}
\Crefname{assumption}{Assumption}{Assumptions}

\AddToHook{env/theorem/before}{\par\addvspace{3pt}}
\AddToHook{env/theorem/after}{\par\addvspace{3pt}}
\AddToHook{env/example/before}{\par\addvspace{3pt}}
\AddToHook{env/example/after}{\par\addvspace{3pt}}
\AddToHook{env/definition/before}{\par\addvspace{3pt}}
\AddToHook{env/definition/after}{\par\addvspace{3pt}}
\AddToHook{env/proposition/before}{\par\addvspace{3pt}}
\AddToHook{env/proposition/after}{\par\addvspace{3pt}}
\AddToHook{env/lemma/before}{\par\addvspace{3pt}}
\AddToHook{env/lemma/after}{\par\addvspace{3pt}}
\AddToHook{env/remark/before}{\par\addvspace{3pt}}
\AddToHook{env/remark/after}{\par\addvspace{3pt}}
\AddToHook{env/corollary/before}{\par\addvspace{3pt}}
\AddToHook{env/corollary/after}{\par\addvspace{3pt}}
\newtheorem{theorem}{Theorem}[section]
\newtheorem{example}[theorem]{Example}
\AddToHook{env/example/begin}{\crefalias{theorem}{example}}
\newtheorem{lemma}[theorem]{Lemma}
\AddToHook{env/lemma/begin}{\crefalias{theorem}{lemma}}
\newtheorem{remark}[theorem]{Remark}
\AddToHook{env/remark/begin}{\crefalias{theorem}{remark}}
\newtheorem{proposition}[theorem]{Proposition}
\AddToHook{env/proposition/begin}{\crefalias{theorem}{proposition}}
\newtheorem{definition}[theorem]{Definition}
\AddToHook{env/definition/begin}{\crefalias{theorem}{definition}}
\newtheorem{conjecture}[theorem]{Conjecture}
\AddToHook{env/conjecture/begin}{\crefalias{theorem}{conjecture}}
\newtheorem{corollary}[theorem]{Corollary}
\AddToHook{env/corollary/begin}{\crefalias{theorem}{corollary}}

\AddToHook{env/question/begin}{\crefalias{theorem}{question}}

\AddToHook{env/assumption/begin}{\crefalias{theorem}{assumption}}

\numberwithin{equation}{section}

\renewcommand{\exp}{\operatorname{exp}}

\usepackage[url=true,backend=biber,style=ext-alphabetic,hyperref=true,giveninits=true]{biblatex}
\title[Tensor products of semistable adelic vector bundles]{On the tensor product of semistable Hermitian adelic vector bundles over an adelic curve}

\usepackage{orcidlink}
\author{Huayi Chen\orcidlink{0000-0001-9199-4386}}
\address{Institute for Theoretical Sciences, Westlake University, No. 600 Dunyu Road, Sandun town, Xihu district, Hangzhou, Zhejiang Province, 310030, China}
\email{chenhuay@westlake.edu.cn}
\urladdr{https://chenhuayi.lab.westlake.edu.cn/}
\author{Yijun Yuan\orcidlink{0000-0001-6571-6980}}
\address{Institute for Theoretical Sciences, Westlake University, No. 600 Dunyu Road, Sandun town, Xihu district, Hangzhou, Zhejiang Province, 310030, China}
\email{941201yuan@gmail.com}
\urladdr{https://yijunyuan.github.io/}

\begin{document}
\subjclass{14G40, 20C30, 14L24, 05E10, 68V20}
\keywords{Arakelov geometry, adelic curves,
Hermitian adelic vector bundles, tensor semistability, Bost conjecture, Lean formalization}
\begin{abstract}
	We prove tensor semistability for Hermitian adelic vector bundles on proper adelic curves whose underlying fields are countable and of characteristic zero. The central estimate concerns the canonical integral Specht lattice of a rectangular partition with fixed height: its logarithmic covolume, divided by its rank and tensor degree, tends to zero. We prove this estimate from monic seminormal norms and an exact hypergeometric law weighted by the square of a Vandermonde product for tableau prefixes. Constant adelic Specht bundles and equivariant contractions then yield the slope inequality for every geometrically semistable subspace of a finite tensor product. The correction is the measure of the Archimedean places times a sublinear Specht term and vanishes in the limit as the auxiliary degree tends to infinity. An instability-filtration argument completes the reduction to arbitrary factors and identifies the tensor product Harder--Narasimhan filtration. The number field specialization gives the Bost tensor semistability conjecture in all ranks.
\end{abstract}
\maketitle

\section{Introduction}\label{sec:introduction}

\subsection{Semistability of tensor vector bundles}
Let $k$ be a field. Assume that $\operatorname{char}(k)=0$, and let $C$ be a regular projective integral curve over $\operatorname{Spec}(k)$. For a nonzero vector bundle $E$ on $C$, define its \emph{slope} by
\[ \operatorname{\mu}(E):=\frac{\operatorname{deg}(c_1(E)\cap[C])}{\operatorname{rk}(E)}, \]
where $c_1(E)$ denotes the first Chern class of $E$, $[C]$ the fundamental cycle of $C$, and $\operatorname{deg}$ the degree of a zero-cycle over $k$. We call $E$ \emph{semistable} if every nonzero vector subbundle $F\subseteq E$ satisfies
\[ \operatorname{\mu}(F)\le\operatorname{\mu}(E). \]
We call $E$ \emph{stable} if the inequality is strict for every nonzero proper subbundle.

A classical result of Narasimhan and Seshadri in complex geometry implies that the tensor product of two semistable vector bundles over a regular projective curve is still semistable. On a compact Riemann surface of genus at least two, their theorem identifies stable bundles of degree zero with irreducible \emph{unitary representations} of the fundamental group. For arbitrary degree, the corresponding statement uses unitary representations of a suitable central extension, or \emph{projectively unitary connections} (cf. \cite[Section 12, Theorem 2]{ns}). Tensoring these connections and passing to \emph{Jordan--H\"older filtrations} gives the assertion for semistable bundles. The tensor product theorem over an arbitrary field of characteristic zero follows by descent to a finitely generated subfield and extension to $\mathbb{C}$; see \cite[Theorem 1.1 and Section 1.1]{bc}.

Ramanan and Ramanathan subsequently established an algebraic extension to \emph{principal bundles} under reductive groups. Over an algebraically closed field of characteristic zero, extension of structure group preserves semistability whenever the homomorphism of connected reductive groups maps the connected centre of the source into the connected centre of the target (cf. \cite[Theorem 3.18]{rr}). Applied to the tensor representation of a product of general linear groups, this recovers the theorem for vector bundles in every genus and all ranks. In characteristic $p>0$, the corresponding unrestricted result requires \emph{strong semistability}, meaning semistability of every iterated Frobenius pullback (cf. \cite[Theorem 3.23]{rr}); tensor products then remain strongly semistable (cf. \cite[Theorem 6.1]{langer}). For ordinary semistable bundles of positive ranks $r$ and $s$, the Ilangovan--Mehta--Parameswaran theorem gives tensor semistability under the bound $r+s\le p+1$ (cf. \cite[Theorem 3.1]{imp}); this rank formulation is recorded in \cite[Section 6, p. 274]{langer}.

\subsection{Arithmetic analogue: Bost's conjecture}

In Arakelov geometry, Stuhler observed that slope semistability can be defined naturally for \emph{Euclidean lattices} (namely the data $\overline{\Lambda}$ of a free abelian group of finite rank $\Lambda$ equipped with a Euclidean inner product on $\Lambda_{\mathbb R}$). If $\overline\Lambda$ is a Euclidean lattice of positive rank $r$, define
\[ \operatorname{\mu}(\overline\Lambda):= -\frac{\ln(\operatorname{covol}(\overline\Lambda))}{r}, \]
where $\operatorname{covol}$ denotes the volume of a fundamental domain (where we consider the Haar measure on $\Lambda_{\mathbb R}$ such that the parallelotope spanned by an orthonormal basis has mass $1$). We call $\overline\Lambda$ semistable if every nonzero sublattice, with its induced metric, has slope at most $\operatorname{\mu}(\overline\Lambda)$. Stuhler also constructed the analogue of the \emph{Harder--Narasimhan filtration} (cf. \cite[Satz 1 and its corollary]{stuhler}). His covolume inequality (cf. \cite[Proposition 2]{stuhler}) supplies the degree inequality underlying this construction; see also \cite[Theorem 1.12]{grayson1}.

Grayson generalized this construction to lattices over rings of integers in number fields, or equivalently to Hermitian vector bundles on their spectra. Every such bundle admits a unique filtration by saturated submodules whose successive quotients, with the quotient metrics, are semistable and have strictly decreasing slopes (cf. \cite[Theorem 1.18 and Definitions 1.19--1.20]{grayson1}). Here a submodule is called \emph{saturated} if its quotient module is torsion-free. Grayson used this filtration in reduction theory for arithmetic groups, and his later work also bounds the maximal slope of a tensor product of two Euclidean lattices by the sum of their maximal slopes plus a constant depending only on their ranks (cf. \cite[Corollary 2.3 and Theorem 5.1]{grayson2}).

The analogy with vector bundles on algebraic curves led Bost to formulate the tensor semistability conjecture. We first recall his slope formalism (cf. \cite[Appendix A.2--A.4]{bost96}). Let $K$ be a number field and let $\mathcal O_K$ be its ring of integers. One defines a \emph{Hermitian vector bundle} over $\operatorname{Spec}(\mathcal O_K)$ to be a finitely generated projective $\mathcal O_K$-module $E$ equipped with conjugation-invariant Hermitian norms on $E\otimes_{\mathcal O_K,\sigma}\mathbb{C}$ for the embeddings $\sigma:K\hookrightarrow\mathbb{C}$, and write $\overline E$ for these data. Assume that $r:=\operatorname{rk}(E)\in\mathbb{Z}_{\ge1}$, and put $M:=\bigwedge_{\mathcal O_K}^r E$. For each embedding $\sigma:K\hookrightarrow\mathbb{C}$, equip $M\otimes_{\mathcal O_K,\sigma}\mathbb{C}$ with the determinant norm $\|\mathord{\cdot}\|_{\sigma,\operatorname{det}}$ induced by the Hermitian norm of $\overline E$. For $s\in M\setminus\{0\}$, define the \emph{normalized Arakelov degree} by
\[ \operatorname{\widehat{deg}}(\overline E):=\frac{1}{[K:\mathbb{Q}]} \left(\ln\bigl(\#(M/\mathcal O_Ks)\bigr) -\sum_{\sigma:K\hookrightarrow\mathbb{C}} \ln(\|s\|_{\sigma,\operatorname{det}})\right). \]
Here $\#(M/\mathcal O_Ks)$ denotes the cardinality of the finite quotient, and the sum runs over all embeddings, including both embeddings in each complex conjugate pair. The product formula shows that this expression is independent of $s$. Define
\[ \operatorname{\widehat{\mu}}(\overline E):=\frac{\operatorname{\widehat{deg}}(\overline E)}{\operatorname{rk}(E)},\qquad \operatorname{\widehat{\mu}}_{\max}(\overline E):= \max_{0\ne F\subseteq E}\operatorname{\widehat{\mu}}(\overline F), \]
where $F$ ranges over nonzero $\mathcal O_K$-submodules with induced metrics. We call $\overline E$ semistable if $\operatorname{\widehat{\mu}}_{\max}(\overline E)=\operatorname{\widehat{\mu}}(\overline E)$. For $K=\mathbb{Q}$, one recovers Stuhler's definition for Euclidean lattices. Bost proposed the following conjecture in his Oberwolfach lecture of 24 July 1997 (cf. \cite{bost97}); see \cite[Problems 1.3--1.5]{bc} for its equivalent formulations.

\begin{conjecture}\label{conj:bost}
	Let $K$ be a number field and let $\mathcal O_K$ be its ring of integers. Assume that $\overline E$ and $\overline F$ are nonzero semistable Hermitian vector bundles on $\operatorname{Spec}(\mathcal O_K)$. Then $\overline E\otimes\overline F$, equipped with the Hermitian tensor product norms, is semistable.
\end{conjecture}

The above conjecture is equivalent to the following equality
\[\operatorname{\widehat{\mu}}_{\max}(\overline E\otimes\overline F)=\operatorname{\widehat{\mu}}_{\max}(\overline E)+\operatorname{\widehat{\mu}}_{\max}(\overline F).\]
Note that the inequality
\[\operatorname{\widehat{\mu}}_{\max}(\overline E\otimes\overline F)\geq\operatorname{\widehat{\mu}}_{\max}(\overline E)+\operatorname{\widehat{\mu}}_{\max}(\overline F)\]
always holds. Hence the conjecture is actually equivalent to the opposite inequality
\[\operatorname{\widehat{\mu}}_{\max}(\overline E\otimes\overline F)\leq\operatorname{\widehat{\mu}}_{\max}(\overline E)+\operatorname{\widehat{\mu}}_{\max}(\overline F)\]

An early general estimate was obtained by Huayi Chen (cf. \cite[Theorem 1.1, formula (2)]{chen09}). For nonzero Hermitian vector bundles $\overline E$ and $\overline F$ over a number field, his theorem gives
\[ \operatorname{\widehat{\mu}}_{\max}(\overline E\otimes\overline F) \le\operatorname{\widehat{\mu}}_{\max}(\overline E)+\operatorname{\widehat{\mu}}_{\max}(\overline F) +\ln(\operatorname{rk}(E))+\ln(\operatorname{rk}(F)). \]
The error depends only on the ranks and is independent of the number field. The proof bounds the degrees of line subbundles by invariant polynomials in the geometrically semistable case, then uses instability filtrations and an explicit form of the Ramanan--Ramanathan theorem to treat the general case (cf. \cite[Sections 4, 6, 8 and 9]{chen09}).

Bost and Chen subsequently sharpened this estimate (cf. \cite[Theorem A]{bc}) to
\[ \operatorname{\widehat{\mu}}_{\max}(\overline E\otimes\overline F) \le\operatorname{\widehat{\mu}}_{\max}(\overline E)+\operatorname{\widehat{\mu}}_{\max}(\overline F)+\frac12\min\{ \ln(\operatorname{rk}(E)),\ln(\operatorname{rk}(F))\}. \]
Their Theorem B (cf. \cite[Theorem B]{bc}) gives
\[ \operatorname{\widehat{\mu}}(\overline V)\le \operatorname{\widehat{\mu}}_{\max}(\overline E)+\operatorname{\widehat{\mu}}_{\max}(\overline F) \]
for every Hermitian subbundle $\overline V\subseteq\overline E\otimes\overline F$ with $\operatorname{rk}(V)\in\{1,2,3,4\}$, and establishes \Cref{conj:bost} for semistable factors with $\operatorname{rk}(E)\operatorname{rk}(F)\le9$ (cf. \cite[Theorem B]{bc}). Their use of the $\varepsilon$-tensor norm, tensor ranks, and geometric semistability also leads to a reduction of the conjecture to subspaces satisfying a geometric-semistability condition (cf. \cite[Problem 1.17, Proposition 7.4 and Theorem 7.5]{bc}).

\subsection{Adelic extensions and known tensor product bounds} \label{subsec:intro-arithmetic-history}
Gaudron extended the arithmetic slope formalism to \emph{pure adelic vector bundles over global fields} (cf. \cite[Sections 3--5]{gaudron}); the purity restriction is made explicit in his errata (cf. \cite{gaudron-errata}). His base field may be a number field or a finite extension of a rational function field over a finite field. The objects are finite-dimensional vector spaces equipped with compatible norms at all places, ultrametric at the non-Archimedean places and standard in a fixed basis outside finitely many places. The Archimedean norms may be general norms. Gaudron constructed the adelic degree and slope polygon, and obtained a canonical Harder--Narasimhan filtration in the Hermitian case (cf. \cite[Proposition 5.13]{gaudron}). For general norms, comparison with the \emph{John and L\"owner Hermitian structures} gives quantitative versions of the Hermitian slope relations (cf. \cite[Sections 4.2 and 5]{gaudron}).

Chen and Moriwaki extended slope theory to \emph{adelic curves}, where the places are parametrized by a measure space and the underlying field need not be global (cf. \cite[Chapters 3--4]{cm}). Let $K$ be a field and let $M_K$ denote the set of absolute values on $K$, including the trivial one. We call a tuple
\[ S:=(K,(\Omega,\mathcal A,\nu),\phi) \]
an adelic curve if $(\Omega,\mathcal A,\nu)$ is a measure space, $\phi:\Omega\to M_K$ is a map, and, on writing $|\mathord{\cdot}|_\omega:=\phi(\omega)$, the function $\omega\mapsto\ln(|a|_\omega)$ is $\mathcal A$-measurable and belongs to $L^1(\Omega,\mathcal A,\nu)$ for every $a\in K^\times:=K\setminus\{0\}$. We call $S$ \emph{proper} if it satisfies the product formula
\[ \forall\,a\in K^\times,\quad \int_\Omega\ln(|a|_\omega)\,\nu(d\omega)=0. \]
Throughout the paper, assume that $S$ is proper and $K$ is countable. For $\omega\in\Omega$, let $K_\omega$ denote the completion of $K$ for $|\mathord{\cdot}|_\omega$. Define
\[ \Omega_\infty:=\{\omega\in\Omega:|\mathord{\cdot}|_\omega \text{ is Archimedean}\}. \]
At every $\omega\in\Omega_\infty$, assume that the absolute value on $K_\omega\simeq\mathbb{R}$ or $\mathbb{C}$ is the usual one. The set $\Omega_\infty$ is $\mathcal A$-measurable and has finite measure (cf. \cite[Propositions 3.1.1--3.1.2]{cm}).

For a finite-dimensional $K$-vector space $E$, put $E_\omega:=E\otimes_K K_\omega$ for $\omega\in\Omega$ and $E^\vee:=\operatorname{Hom}_K(E,K)$. We define an \emph{adelic vector bundle} on $S$ to be a pair $\overline E:=(E,\xi)$, where $\xi=(\|\mathord{\cdot}\|_\omega)_{\omega\in\Omega}$ is a family of norms on the spaces $E_\omega$ satisfying the following measurability and domination conditions. For every $s\in E\setminus\{0\}$ and every $\ell\in E^\vee\setminus\{0\}$, the functions
\[ \omega\longmapsto\ln(\|s\|_\omega),\qquad \omega\longmapsto\ln(\|\ell\|_{\omega,*}) \]
are $\mathcal A$-measurable and each has an upper bound in $L^1(\Omega,\mathcal A,\nu)$ for $\nu$-almost every $\omega$; here $\|\mathord{\cdot}\|_{\omega,*}$ denotes the dual norm. We say that $\overline E$ is \emph{Hermitian} if its norms are ultrametric at non-Archimedean places and induced by Euclidean or Hermitian inner products at Archimedean places. These definitions and the equivalent basis criterion for domination are recalled in \Cref{sec:slopes}; see \cite[Definitions 4.1.1, 4.1.2 and 4.1.28]{cm}. The Hermitian Harder--Narasimhan filtration is constructed in \cite[Propositions 4.3.37--4.3.38]{cm}.

Assume that $\operatorname{char}(K)=0$. Let $N\in\mathbb{Z}_{\ge2}$ and let $(\overline E_j)_{j=1}^N$ be positive-rank Hermitian adelic vector bundles on $S$. Equip their tensor product with the \emph{orthogonal tensor product}: \emph{Hilbert--Schmidt tensor norms} at Archimedean places and \emph{$\varepsilon$-tensor norms} at the other places. Define the \emph{minimal slope} of a positive-rank adelic bundle to be the infimum of the slopes of its nonzero quotients equipped with quotient norms. Chen and Moriwaki proved (cf. \cite[Corollary 5.6.2, formula (5.23)]{cm})
\begin{equation}\label{eq:intro-cm-tensor-bound}
	\operatorname{\widehat{\mu}}_{\min}\left(\bigotimes_{j=1}^N\overline E_j\right) \ge\sum_{j=1}^N\operatorname{\widehat{\mu}}_{\min}(\overline E_j) -\frac12\nu(\Omega_\infty)\sum_{j=1}^N\ln(\operatorname{rk}(E_j)).
\end{equation}
For general adelic bundles with their $\varepsilon,\pi$-tensor norms, formula (5.21) of the same corollary has coefficient $\tfrac32\nu(\Omega_\infty)$. The proof proceeds in two stages. First, consider a nonzero quotient of the tensor product whose determinant point is geometrically semistable. For such a quotient, \cite[Proposition 5.3.1, formula (5.7)]{cm} bounds its slope below by the sum of the slopes of the factors, minus an error term. For the orthogonal Hermitian tensor product, \cite[Remark 5.3.2, formula (5.9)]{cm} gives this error as
\[ \frac12\nu(\Omega_\infty) \sum_{j=1}^N\ln(\operatorname{rk}(E_j)). \]
The proof then uses instability filtrations to treat arbitrary nonzero quotients. The resulting comparison (cf. \cite[Theorem 5.6.1, formula (5.16)]{cm}) bounds the slope of each quotient below by the average weight of the quotient filtration induced by the tensor product of the factors' Harder--Narasimhan filtrations, minus the same error. Here weights are counted with their graded dimensions. The smallest weight of each factor's Harder--Narasimhan filtration is its minimal slope. Hence every weight in the induced quotient filtration is at least the sum of these minimal slopes. Taking the infimum over all nonzero quotients gives \eqref{eq:intro-cm-tensor-bound}.

\subsection{The main statement}
In this paper, we prove that the Archimedean error can be removed in the Hermitian case. The main statement is the following adelic extension of \Cref{conj:bost}.

\begin{theorem}\label{thm:bost}
	Let $S:=(K,(\Omega,\mathcal A,\nu),\phi)$ be a proper adelic curve. Assume that $K$ is countable and $\operatorname{char}(K)=0$. Let $\overline E$ and $\overline F$ be nonzero Hermitian adelic vector bundles on $S$. Then
	\begin{equation}\label{eq:main}
		\operatorname{\widehat{\mu}}_{\max}(\overline E\otimes\overline F) =\operatorname{\widehat{\mu}}_{\max}(\overline E) +\operatorname{\widehat{\mu}}_{\max}(\overline F).
	\end{equation}
	Consequently, if $\overline E$ and $\overline F$ are semistable, then $\overline E\otimes\overline F$ is semistable. The same conclusions hold for every nonempty finite family of nonzero Hermitian adelic vector bundles on $S$.
\end{theorem}

The key to the proof is an asymptotic slope estimate for auxiliary adelic bundles constructed from rectangular integral Specht lattices.

For $(a,k)\in\mathbb{Z}_{\ge1}^2$, put $d:=ak$ and $\lambda:=(k,\ldots,k)\in\mathbb N^a$. Let $L_{a,d}$ be the full \emph{integral Specht lattice} of shape $\lambda$, realized by unnormalized column anti-symmetrization in $(\mathbb{Z}^a)^{\otimes d}$ and equipped with the induced Euclidean metric. Its integral structure is the orbit lattice of all tableaux, as defined in \Cref{sec:specht}. By extension of scalars, define its associated Hermitian adelic vector bundle $\overline L_{a,d}^{\,S}$ with underlying $K$-vector space
\[ L_{a,d}^{\,S}:=K\otimes_\mathbb{Z} L_{a,d}. \]
At real Archimedean places retain the Euclidean metric, at complex Archimedean places take its Hermitian extension, and at non-Archimedean places use the supremum norm in any $\mathbb{Z}$-basis of $L_{a,d}$ (see \Cref{subsec:finite-groups}). By \eqref{eq:lattice-extension-slope},
\[ -\operatorname{\widehat{\mu}}(\overline L_{a,d}^{\,S}) =\nu(\Omega_\infty)\frac{\ln(\operatorname{covol}(L_{a,d}))}{\#(\operatorname{SYT}(\lambda))}. \]
For fixed $S$ and $a$, the estimate used to remove the tensor-slope error is
\begin{equation}\label{eq:intro-sublinear}
	\lim\limits_{k\to\infty} \frac{-\operatorname{\widehat{\mu}}(\overline L_{a,ak}^{\,S})}{ak}=0,
\end{equation}
proved in \Cref{prop:sublinear}.

The standard Specht vectors and Young's seminormal representation are classical (cf. \cite[Section 2.1]{kn}, \cite[Section 2, equation (1)]{rh}). Their normalization matters here: multiplying each basis vector by a positive scalar changes the covolume by that scalar raised to the rank. We therefore check the initial seminormal norm in the tensor metric. The standard vectors generate a sublattice of finite index in $L_{a,d}$; their Gram matrix gives an upper bound for its covolume. The additional estimate is obtained by averaging the logarithms of the squared seminormal norms. The prefix shape of a uniformly chosen standard tableau has an exact hypergeometric distribution weighted by a squared Vandermonde determinant. This weight suppresses small gaps between shifted row lengths and makes the average norm cost per tensor position tend to zero. \Cref{sec:seminormal,sec:asymptotic} prove the probability law, its normalization, the uniform gap bound, and the resulting covolume estimate.

The arithmetic step uses $\overline L_{a,d}^{\,S}$. Absolute irreducibility makes the auxiliary bundles semistable. The contraction maps preserve the integral norm bounds at non-Archimedean places. At Archimedean places, we retain the determinant-tensor factors of \cite[Remark A.2.9]{cmp}, so their contribution is measured by the Specht covolumes. Assume that $\operatorname{char}(K)=0$, let $N\in\mathbb{Z}_{\ge1}$, and let $(\overline E_j)_{j=1}^N$ be positive-rank Hermitian adelic bundles on $S$. If a nonzero subspace $V\subseteq\bigotimes_{j=1}^N E_j$ has geometrically semistable Pl\"ucker point for $\prod_{j=1}^N\operatorname{SL}(E_j)$, the slope comparison at a fixed auxiliary degree has correction
\[ -\frac1d\left( \operatorname{\widehat{\mu}}(\overline L_{\operatorname{dim}_K(V),d}^{\,S}) +\sum_{j=1}^N\operatorname{\widehat{\mu}}(\overline L_{\operatorname{dim}_K(E_j),d}^{\,S}) \right). \]
Here the auxiliary degree $d$ ranges over sufficiently divisible positive integers. With all ranks fixed, \eqref{eq:intro-sublinear} makes this correction tend to zero. \Cref{thm:arithmetic-git} consequently gives
\begin{equation}\label{eq:intro-git}
	\operatorname{\widehat{\mu}}(\overline V)\le \sum_{j=1}^N\operatorname{\widehat{\mu}}(\overline E_j).
\end{equation}
This is the geometric-semistability inequality formulated in \cite[Problem 1.17]{bc}.

To pass from \eqref{eq:intro-git} to arbitrary tensor subspaces, we use the instability filtrations and rank induction of \cite[Theorems 5.4.3 and 5.6.1]{cm}. \Cref{sec:allrank} gives the required first variation comparison on graded quotients and completes the simultaneous induction for semistability and the Harder--Narasimhan filtration. In the number field case, \cite[Proposition 7.4 and Theorem 7.5]{bc} provide another reduction from \eqref{eq:intro-git} to the tensor product conjecture.

\subsection{Lean 4 formalization}
The results of this paper have been formalized in Lean 4 (cf. \cite{mouraLean4TheoremProver2021}, \href{https://github.com/leanprover/lean4/releases/tag/v4.35.0}{version 4.35.0}) using the mathlib4 library (cf. \cite{mathlib2020}, \href{https://github.com/leanprover-community/mathlib4/commit/56a005855ea17e0b0958f73ec801ae140b0a1679}{commit \texttt{56a0058}}). The code is available at \url{https://github.com/YijunYuan/AdelicBostConjecture/}. The formalization is \texttt{sorry}-free and depends only on the standard Lean axioms \texttt{Classical.choice}, \texttt{propext}, and \texttt{Quot.sound}.
\begin{remark}
	The formalized statements faithfully reflect the mathematical statements in this paper, except for a few intermediate results whose statements have been weakened to simplify the code. We refer the reader to \url{https://github.com/YijunYuan/AdelicBostConjecture#statement-coverage} for a detailed discussion of statement coverage. These weakened statements do not affect the main results of this paper.
\end{remark}
This formalization depends on several external libraries:
\begin{enumerate}
	\item The foundations of Arakelov geometry from \cite{cm} used in this work have been formalized by Jiedong Jiang and Yijun Yuan in the separate \href{https://github.com/YijunYuan/ArakelovGeometry/tree/lite}{ArakelovGeometry} library.
	\item Harder--Narasimhan theory has been formalized by Yijun Yuan in the separate \href{https://github.com/YijunYuan/HarderNarasimhan}{HarderNarasimhan} library (cf. \cite{yuan2026formalizationhardernarasimhantheory}). This formalization is based on Huayi Chen and Marion Jeannin's paper \cite{chen2023hardernarasimhangames} and replaces the basic construction of Harder--Narasimhan theory in \cite{cm}.
	\item We adapt selected definitions and proofs from the following formalization projects and mathlib contributions:
	\begin{enumerate}
		\item \href{https://github.com/openai/math}{\texttt{openai/math}}: tableau and word models, polytabloids, and the irreducibility of complex Specht modules from its Saxl development.
		\item \href{https://github.com/TauCetiProject/TauCeti}{TauCeti}: lattice covolumes and Gram determinants, Garnir relations and rational straightening for Specht modules, Schur--Weyl duality, and symmetric-power polarization.
		\item The \href{https://github.com/mathlib-initiative/Etingof-RepresentationTheory-draft1}{Etingof representation theory project}: partition and Young diagram constructions, and the hook-length formula for standard Young tableaux.
		\item \href{https://github.com/AxiomMath/Stacks}{\texttt{AxiomMath/Stacks}}: the determinant covariance of alternating tensors under tensor powers of matrices.
		\item \href{https://github.com/leanprover-community/mathlib4/pull/40948}{Mathlib PR \#40948} by Nailin Guan: base change for exterior powers and the projection from the exterior algebra onto a homogeneous exterior power.
		\item \href{https://github.com/leanprover-community/mathlib4/pull/44845}{Mathlib PR \#44845} by Daniel Horton: tensor products of finite families of projective modules, together with tensor--Hom and duality equivalences for finite projective modules.
	\end{enumerate}
	The source revisions, authorship and license notices, and local adaptations are recorded in the repository's source files and provenance records.
\end{enumerate}

The formalization was initially carried out using \href{https://github.com/frenzymath/Archon-Horizon}{\texttt{Archon-Horizon}}, an agentic auto-formalization system developed by the AI4Math team at BICMR, Peking University. We used \texttt{GPT-6 Astra} (cf. \cite{GPT6Astra}) from OpenAI as the base model. The authors subsequently refactored and polished the code extensively.

\subsection{AI usage declaration}
\texttt{GPT-6 Astra} (cf. \cite{GPT6Astra}) from OpenAI was used to assist with mathematical discussions and the development of arguments, drafting and revising the exposition, checking bibliographic references, editing the \LaTeX{} source and formalizing the proofs. The authors retain full responsibility for the mathematical claims, proofs, and final text.

\section{Preliminaries on Hermitian adelic vector bundles} \label{sec:slopes}

\subsection{Proper adelic curves}\label{subsec:adelic-curves}

Let $K$ be a field. Assume that $K$ is countable. Define $M_K$ to be the set of absolute values on $K$, including the trivial absolute value. The following definitions are those of \cite[Section 3.1]{cm}.

\begin{definition}
	Let $(\Omega,\mathcal A,\nu)$ be a measure space and let $\phi:\Omega\to M_K$ be a map. For $\omega\in\Omega$, define $|\mathord{\cdot}|_\omega:=\phi(\omega)$. We call
	\[ S:=(K,(\Omega,\mathcal A,\nu),\phi) \]
	an \emph{adelic curve} if, for every $a\in K^\times:=K\setminus\{0\}$, the function $\omega\mapsto\ln(|a|_\omega)$ is $\mathcal A$-measurable and belongs to $L^1(\Omega,\mathcal A,\nu)$. We call $S$ \emph{proper} if it satisfies the \emph{product formula}
	\begin{equation}\label{eq:product-formula}
		\forall\,a\in K^\times,\quad \int_\Omega\ln(|a|_\omega)\,\nu(d\omega)=0.
	\end{equation}
\end{definition}

Throughout the paper, assume that $S$ is proper and $K$ is countable. For $\omega\in\Omega$, let $K_\omega$ be the completion of $K$ with respect to $|\mathord{\cdot}|_\omega$. Define
\[ \Omega_\infty:=\{\omega\in\Omega:|\mathord{\cdot}|_\omega \text{ is Archimedean}\}. \]
Following the normalization of \cite[Section 4.1]{cm}, assume that, for every $\omega\in\Omega_\infty$, the absolute value on $K_\omega\simeq\mathbb{R}$ or $\mathbb{C}$ is the usual one. The set $\Omega_\infty$ is $\mathcal A$-measurable and has finite measure (cf. \cite[Propositions 3.1.1--3.1.2]{cm}). If $\operatorname{char}(K)>0$, then $\Omega_\infty=\varnothing$ (cf. \cite[Section 3.1]{cm}).

\subsection{Hermitian norm families and adelic vector bundles} \label{subsec:norm-families}

Let $E$ be a $K$-vector space of dimension $r\in\mathbb{Z}_{\ge1}$. For $\omega\in\Omega$, define $E_\omega:=E\otimes_K K_\omega$. We identify $E$ with its image in $E_\omega$ under $s\mapsto s\otimes1$. We define a \emph{norm family} on $E$ to be a family $\xi=(\|\mathord{\cdot}\|_\omega)_{\omega\in\Omega}$ such that, for every $\omega\in\Omega$, the map $\|\mathord{\cdot}\|_\omega:E_\omega\to\mathbb{R}_{\ge0}$ is a norm on the $K_\omega$-vector space $E_\omega$ with respect to $|\mathord{\cdot}|_\omega$. We call $\xi$ \emph{Hermitian} if $\|\mathord{\cdot}\|_\omega$ is \emph{ultrametric} for $\omega\in\Omega\setminus\Omega_\infty$ and is induced by a positive-definite real or Hermitian inner product for $\omega\in\Omega_\infty$. Here ultrametric means
\[ \forall\,(x,y)\in E_\omega\times E_\omega,\quad \|x+y\|_\omega\le\max\{\|x\|_\omega,\|y\|_\omega\}. \]
Write $E^\vee:=\operatorname{Hom}_K(E,K)$. Define the \emph{dual norm family} $\xi^\vee=(\|\mathord{\cdot}\|_{\omega,*})_{\omega\in\Omega}$ by
\[ \forall\,\ell\in(E^\vee)_\omega,\quad \|\ell\|_{\omega,*}:= \sup_{x\in E_\omega\setminus\{0\}} \frac{|\ell(x)|_\omega}{\|x\|_\omega}. \]

\begin{definition}
	We call a norm family $\xi=(\|\mathord{\cdot}\|_\omega)_{\omega\in\Omega}$ on $E$ \emph{measurable} if $\omega\mapsto\|s\|_\omega$ is $\mathcal A$-measurable for every $s\in E$. We call $\xi$ \emph{upper dominated} if, for every $s\in E\setminus\{0\}$, there exists a function $g_s:\Omega\to\mathbb{R}$ with $g_s\in L^1(\Omega,\mathcal A,\nu)$ such that
	\[ \ln(\|s\|_\omega)\le g_s(\omega) \quad\text{for $\nu$-almost every $\omega$}. \]
	We call $\xi$ \emph{dominated} if both $\xi$ and $\xi^\vee$ are upper dominated. We call a pair $\overline E:=(E,\xi)$ a \emph{Hermitian adelic vector bundle} on $S$ if $\xi$ is Hermitian and dominated, and both $\xi$ and $\xi^\vee$ are measurable. Define its \emph{rank} by $\operatorname{rk}(\overline E):=\operatorname{dim}_K(E)$.
\end{definition}

These are the conventions of \cite[Definitions 4.1.1, 4.1.2 and 4.1.28]{cm}. A useful equivalent form of domination is obtained by comparison with a basis. Let $\mathbf e=(e_1,\ldots,e_r)$ be a $K$-basis of $E$. Define the \emph{reference norm family} $\xi_{\mathbf e}=(\|\mathord{\cdot}\|_{\mathbf e,\omega})_{\omega\in\Omega}$ by
\[
	\begin{gathered}
		\forall\,(x_1,\ldots,x_r)\in K_\omega^r,\\
		\left\|\sum_{i=1}^r x_i e_i\right\|_{\mathbf e,\omega}:=
		\begin{cases}
			\max_{i\in\{1,\ldots,r\}}|x_i|_\omega,
			 & \omega\notin\Omega_\infty, \\[2pt]
			\bigl(\sum_{i=1}^r|x_i|_\omega^2\bigr)^{1/2},
			 & \omega\in\Omega_\infty,
		\end{cases}
	\end{gathered}
\]
For a measurable Hermitian norm family $\xi$, define
\[ d_\omega(\xi,\xi_{\mathbf e}):= \sup_{x\in E_\omega\setminus\{0\}} \bigl|\ln(\|x\|_\omega)-\ln(\|x\|_{\mathbf e,\omega})\bigr|. \]
The family $\xi$ is dominated if and only if
\begin{equation}\label{eq:domination-distance}
	\int_\Omega d_\omega(\xi,\xi_{\mathbf e})\,\nu(d\omega)<\infty
\end{equation}
for one, or equivalently, for every $K$-basis $\mathbf e$ (cf. \cite[Corollary 4.1.10 and Remark 4.1.12]{cm}). The Hermitian norms agree with their biduals. Moreover, changing between the usual coordinate reference norms at infinity changes the local distance by a bounded amount supported on $\Omega_\infty$, whose measure is finite. Thus the criterion applies to the displayed Euclidean reference norm. The distance function is measurable: since $K$ is countable and dense in each $K_\omega$, its defining supremum can be taken over the countable set $E\setminus\{0\}$. The same observation shows that a measurable Hermitian norm family has measurable dual, since, for $\ell\in E^\vee$,
\[ \|\ell\|_{\omega,*} =\sup_{s\in E\setminus\{0\}} \frac{|\ell(s)|_\omega}{\|s\|_\omega}. \]

\subsection{Subquotients, duals, and tensor products} \label{subsec:bundle-operations}

Let $E'\subseteq E$ be a $K$-vector subspace and put $Q:=E/E'$. For $\omega\in\Omega$, equip $E'_\omega:=E'\otimes_K K_\omega$ with the restriction of $\|\mathord{\cdot}\|_\omega$. For $q\in Q_\omega$, define its \emph{quotient norm} by
\[ \|q\|_{Q,\omega}:=\inf_{x\in q} \|x\|_\omega. \]
On vector subspaces we always consider induced norms; on quotient vector spaces we always consider quotient norms. If $\overline F$ is another Hermitian adelic vector bundle, we define their orthogonal tensor product by equipping $E\otimes_K F$ with the Hilbert--Schmidt tensor norm at each Archimedean place and the $\varepsilon$-tensor norm at each non-Archimedean place. The Hilbert--Schmidt tensor norm is also called the orthogonal tensor product norm in \cite[Section 1.2.9]{cm}. Explicitly, at a non-Archimedean place, define the $\varepsilon$-tensor norm by
\[ \forall\,z\in E_\omega\otimes_{K_\omega}F_\omega,\quad \|z\|_{\omega,\varepsilon}:= \sup_{(\ell,m)\in (E_\omega^\vee\setminus\{0\})\times(F_\omega^\vee\setminus\{0\})} \frac{|(\ell\otimes m)(z)|_\omega} {\|\ell\|_{\omega,*}\|m\|_{\omega,*}}. \]
At an Archimedean place, define its inner product by
\[ \forall\,(x,x',y,y')\in E_\omega^2\times F_\omega^2,\quad \langle x\otimes y,x'\otimes y'\rangle_\omega :=\langle x,x'\rangle_{E,\omega} \langle y,y'\rangle_{F,\omega}. \]

Define $\operatorname{det}(E):=\bigwedge\nolimits^r E$ and give it the \emph{determinant norms} of \cite[Section 1.1.13]{cm}. At an Archimedean place they satisfy
\[ \forall\,(x_1,\ldots,x_r)\in E_\omega^r,\quad \|x_1\wedge\cdots\wedge x_r\|_{\omega,\operatorname{det}}^2 =\operatorname{det}\bigl((\langle x_i,x_j\rangle_\omega)_{i,j=1}^r\bigr). \]
At a non-Archimedean place, for $z\in\operatorname{det}(E_\omega)\setminus\{0\}$, define them by
\[ \|z\|_{\omega,\operatorname{det}}:= \inf_{\substack{(x_1,\ldots,x_r)\in E_\omega^r\\ x_1\wedge\cdots\wedge x_r=z}} \prod_{i=1}^r\|x_i\|_\omega. \]
Write $\operatorname{det}(\overline E)$ for the resulting adelic line bundle.

\begin{proposition}\label{prop:adelic-operations}
	Subspaces, quotients, duals, orthogonal tensor products, and determinant lines of Hermitian adelic vector bundles on $S$ are Hermitian adelic vector bundles. For a $K$-vector subspace $F\subseteq E$, the canonical identifications
	\[ (E/F)^\vee\simeq F^\perp\subseteq E^\vee, \qquad F^\vee\simeq E^\vee/F^\perp, \qquad F^\perp:=\{\ell\in E^\vee:\ell|_F=0\}, \]
	are isometries for the induced and quotient norms.
\end{proposition}
\begin{proof}
	Domination is preserved by these operations by \cite[Proposition 4.1.19]{cm}. Their measurability follows from \cite[Proposition 4.1.24 and Remark 4.1.34]{cm}, since the countable field $K$ is dense in every $K_\omega$. Thus \cite[Proposition 4.1.32]{cm} applies. The local norms remain ultrametric or induced by inner products, respectively. The dual identifications are the local subspace and quotient duality isometries; see \cite[Proposition 1.1.20 and Section 1.2.3]{cm}.
\end{proof}

\subsection{Arakelov degree and elementary slope comparisons} \label{subsec:degrees}

For an adelic line bundle $\overline L$ and $s\in L\setminus\{0\}$, define
\[ \operatorname{\widehat{deg}}(\overline L):= -\int_\Omega\ln(\|s\|_{L,\omega})\,\nu(d\omega). \]
Domination of the norm family and its dual makes the integrand integrable. If $s'=as$ with $a\in K^\times$, then \eqref{eq:product-formula} gives
\[ -\int_\Omega\ln(\|s'\|_{L,\omega})\,\nu(d\omega) =-\int_\Omega\ln(\|s\|_{L,\omega})\,\nu(d\omega) -\int_\Omega\ln(|a|_\omega)\,\nu(d\omega) =\operatorname{\widehat{deg}}(\overline L). \]
For a positive-rank Hermitian adelic vector bundle $\overline E$, define its \emph{Arakelov degree} and \emph{slope} by
\[ \operatorname{\widehat{deg}}(\overline E):=\operatorname{\widehat{deg}}(\operatorname{det}(\overline E)), \qquad \operatorname{\widehat{\mu}}(\overline E):= \frac{\operatorname{\widehat{deg}}(\overline E)}{\operatorname{rk}(\overline E)}. \]
Set $\operatorname{\widehat{deg}}(0):=0$.

We recall the following degree formulas and contraction comparison.

Let $\overline E$ and $\overline F$ be positive-rank Hermitian adelic vector bundles. The duality formula (cf. \cite[Proposition 4.3.10, equation (4.25)]{cm}) gives
\begin{equation}\label{eq:degree-duality}
	\operatorname{\widehat{deg}}(\overline E^\vee)=-\operatorname{\widehat{deg}}(\overline E).
\end{equation}
The tensor product degree formula (cf. \cite[Proposition 4.3.19, equation (4.35)]{cm}) and its slope form are
\begin{align}
	\operatorname{\widehat{deg}}(\overline E\otimes\overline F)
	 & =\operatorname{rk}(\overline F)\operatorname{\widehat{deg}}(\overline E)
	+\operatorname{rk}(\overline E)\operatorname{\widehat{deg}}(\overline F),
	\label{eq:tensor-degree}                                                                 \\
	\operatorname{\widehat{\mu}}(\overline E\otimes\overline F)
	 & =\operatorname{\widehat{\mu}}(\overline E)+\operatorname{\widehat{\mu}}(\overline F).
	\label{eq:tensor-slope}
\end{align}
For a subspace $A\subseteq E$ with induced norms and quotient $Q:=E/A$ with quotient norms, degree additivity (cf. \cite[Proposition 4.3.13, equation (4.29)]{cm}) gives
\begin{equation}\label{eq:degree-additivity}
	\operatorname{\widehat{deg}}(\overline E) =\operatorname{\widehat{deg}}(\overline A)+\operatorname{\widehat{deg}}(\overline Q).
\end{equation}
If $u:E\to F$ is injective and $\|u\|_\omega\le1$ for every $\omega\in\Omega$, its image $J:=u(E)$ with induced norms satisfies the contraction comparison (cf. \cite[Proposition 4.3.18, inequality (4.32)]{cm}):
\begin{equation}\label{eq:contract-degree}
	\operatorname{\widehat{\mu}}(\overline E)\le\operatorname{\widehat{\mu}}(\overline J).
\end{equation}
Here $\|u\|_\omega$ denotes the \emph{operator norm} of the scalar extension $u_\omega:E_\omega\to F_\omega$.

For a positive-rank Hermitian adelic vector bundle $\overline E$, define
\[ \Theta(E):=\{F\subseteq E:F\text{ is a $K$-vector subspace of }E\}. \]
Define its \emph{maximal slope} and \emph{minimal slope} by
\[ \operatorname{\widehat{\mu}}_{\max}(\overline E):= \sup_{\substack{F\in\Theta(E)\\F\ne0}}\operatorname{\widehat{\mu}}(\overline F),\qquad \operatorname{\widehat{\mu}}_{\min}(\overline E):= \inf_{\substack{F\in\Theta(E)\\F\ne E}} \operatorname{\widehat{\mu}}(\overline E/\overline F). \]
The subspaces have induced norms and the quotients have quotient norms. These are finite real numbers, and both extrema are attained (cf. \cite[Proposition 4.3.28 and Corollary 4.3.27]{cm}). We call $\overline E$ \emph{semistable} if $\operatorname{\widehat{\mu}}_{\max}(\overline E)=\operatorname{\widehat{\mu}}(\overline E)$. By \eqref{eq:degree-additivity}, this is equivalent to the following condition: for every nonzero quotient $q:E\twoheadrightarrow Q$, one has
\[ \operatorname{\widehat{\mu}}(\overline Q)\ge\operatorname{\widehat{\mu}}(\overline E). \]
Duality gives
\begin{equation}\label{eq:slope-duality}
	\operatorname{\widehat{\mu}}_{\min}(\overline E) =-\operatorname{\widehat{\mu}}_{\max}(\overline E^\vee).
\end{equation}
Together with the quotient criterion, this shows that the dual of a semistable bundle is semistable. For a semistable bundle both extreme slopes equal its slope.

Let $u:E\to F$ be a nonzero $K$-linear map between positive-rank Hermitian adelic vector bundles. The function $\omega\mapsto\ln(\|u\|_\omega)$ is measurable and integrable by \cite[Remark 4.1.12 and Proposition 4.3.16]{cm}. Define its \emph{height} by
\[ h(u):=\int_\Omega\ln(\|u\|_\omega)\,\nu(d\omega). \]

\begin{proposition}\label{prop:slope-comparison}
	Let $u:E\to F$ be a nonzero $K$-linear map between positive-rank Hermitian adelic vector bundles. Then
	\[ \operatorname{\widehat{\mu}}_{\min}(\overline E) \le\operatorname{\widehat{\mu}}_{\max}(\overline F)+h(u). \]
	If $u$ is injective, one has
	\[ \operatorname{\widehat{\mu}}_{\max}(\overline E) \le\operatorname{\widehat{\mu}}_{\max}(\overline F)+h(u). \]
	If $u$ is surjective, one has
	\[ \operatorname{\widehat{\mu}}_{\min}(\overline E) \le\operatorname{\widehat{\mu}}_{\min}(\overline F)+h(u). \]
\end{proposition}
\begin{proof}
	Apply \cite[Proposition 4.3.31]{cm} to $u:E\to F$. Its assertion (3) gives the inequality for a nonzero map, assertion (1) gives the injective case, and assertion (2) gives the surjective case.
\end{proof}

\subsection{The Harder--Narasimhan filtration} \label{subsec:hn}

We recall the Harder--Narasimhan theorem in the Hermitian setting.

\begin{theorem}\label{thm:hn-background}
	For every positive-rank Hermitian adelic vector bundle $\overline E$, there exist a unique integer $\ell\in\mathbb{Z}_{\ge1}$ and a unique flag
	\[ 0=E_0\subsetneq E_1\subsetneq\cdots \subsetneq E_\ell=E \]
	such that, for each $i\in\{1,\ldots,\ell\}$, the bundle $\overline G_i:=\overline E_i/\overline E_{i-1}$ is semistable and, with $\mu_i:=\operatorname{\widehat{\mu}}(\overline G_i)$,
	\[ \mu_1>\cdots>\mu_\ell. \]
	Moreover,
	\[ \operatorname{\widehat{\mu}}_{\max}(\overline E)=\mu_1,\qquad \operatorname{\widehat{\mu}}_{\min}(\overline E)=\mu_\ell,\qquad \operatorname{\widehat{deg}}(\overline E)=\sum_{i=1}^{\ell}\operatorname{dim}_K(G_i)\mu_i. \]
\end{theorem}
\begin{proof}
	For $\overline E$, the construction in \cite[Section 4.3.10 and Proposition 4.3.37]{cm} gives a flag with semistable grades of strictly decreasing slopes. Uniqueness follows from \cite[Proposition 4.3.38]{cm}. In that construction, $E_1$ is the greatest subspace of maximal slope, so $\mu_1=\operatorname{\widehat{\mu}}_{\max}(\overline E)$. The identification $\mu_\ell=\operatorname{\widehat{\mu}}_{\min}(\overline E)$ follows from \cite[Proposition 4.3.41]{cm}. The degree identity follows from \eqref{eq:degree-additivity} applied to the successive terms of the flag.
\end{proof}

We call this flag the \emph{Harder--Narasimhan flag} of $\overline E$. Define its associated decreasing, left-continuous real filtration by
\[ \forall\,t\in\mathbb{R},\quad \mathcal H_E^t E:=\sum_{\substack{i\in\{1,\ldots,\ell\}\\\mu_i\ge t}}E_i. \]
In particular $\mathcal H_E^t E=E$ for $t\le\mu_\ell$ and $\mathcal H_E^t E=0$ for $t>\mu_1$. Uniqueness shows that every isometry preserves this filtration. More generally, \cite[Proposition 4.3.49]{cm} gives, for a nonzero map $u:E\to F$,
\[ \forall\,t\in\mathbb{R},\quad u(\mathcal H_E^t E)\subseteq\mathcal H_F^{\,t-h(u)}F. \]
For a contraction at every place, $h(u)\le0$, so in particular $u(\mathcal H_E^t E)\subseteq\mathcal H_F^t F$.

\subsection{Finite group actions and lattice specializations} \label{subsec:finite-groups}

\begin{proposition}\label{prop:finitegroup} \label{prop:arithmetic-finitegroup}
	Let a finite group $G$ act $K$-linearly on a positive-rank Hermitian adelic vector bundle $\overline E$, and suppose that
	\[ \forall\,\omega\in\Omega,\quad \forall\,(g,x)\in G\times E_\omega,\quad \|gx\|_{E,\omega}=\|x\|_{E,\omega}. \]
	If $E$ is irreducible as a representation of $G$ over $K$, then $\overline E$ is semistable. Let $N\in\mathbb{Z}_{\ge1}$ and let $(G_j,\overline E_j)_{j=1}^N$ be such actions with $E_j$ absolutely irreducible for every $j\in\{1,\ldots,N\}$. Then $\bigotimes_{j=1}^N\overline E_j$ and its dual are semistable for the action of $\prod_{j=1}^N G_j$.
\end{proposition}
\begin{proof}
	By the uniqueness of the Harder--Narasimhan flag in \Cref{thm:hn-background}, which follows from \cite[Proposition 4.3.38]{cm}, every element of $G$ preserves the Harder--Narasimhan flag of $\overline E$. If $\overline E$ is unstable, its first nonzero term is a proper $G$-invariant subspace, contradicting irreducibility. For the second assertion, each $E_j$ becomes irreducible over an algebraic closure of $K$. Applying \cite[Proposition 2.3.23, pp. 65 and 107]{kowalski} successively over that field shows that the external tensor product $\bigotimes_{j=1}^N E_j$ is absolutely irreducible for $\prod_{j=1}^N G_j$. Its dual is also absolutely irreducible: over the algebraic closure, the annihilator of a nonzero proper invariant subspace of the dual would be a nonzero proper invariant subspace of the tensor product. The local tensor norms and dual norms are preserved by the product action. Applying the first assertion to these two representations proves the claim.
\end{proof}

We call a finite free $\mathbb{Z}$-module $L$ equipped with a positive-definite inner product on $L\otimes_\mathbb{Z}\mathbb{R}$ a Euclidean lattice. Suppose that $r:=\operatorname{rk}(L)\in\mathbb{Z}_{\ge1}$. Write $\operatorname{covol}(L)$ for the volume of a fundamental parallelepiped in $L\otimes_\mathbb{Z}\mathbb{R}$, and define its Euclidean degree and slope by
\[ \operatorname{deg}(L):=-\ln(\operatorname{covol}(L)),\qquad \operatorname{\mu}(L):=\frac{\operatorname{deg}(L)}r. \]

We define a Hermitian adelic vector bundle $\overline L^{\,S}$ on the $K$-vector space $L^{\,S}:=K\otimes_\mathbb{Z} L$ by extension of scalars as follows. At an Archimedean place $\omega$, extend the given inner product to a real or Hermitian inner product on $K_\omega\otimes_\mathbb{Z} L$. At a non-Archimedean place $\omega$, choose any $\mathbb{Z}$-basis $(e_1,\ldots,e_r)$ of $L$ and define
\begin{equation}\label{eq:scalar-extension-lattice-norm}
	\forall\,(x_1,\ldots,x_r)\in K_\omega^r,\quad \left\|\sum_{i=1}^r x_i(1\otimes e_i)\right\|_{L,\omega} :=\max_{i\in\{1,\ldots,r\}}|x_i|_\omega.
\end{equation}
An integral change of basis and its inverse have entries of absolute value at most one at every non-Archimedean place. Hence this norm is independent of the chosen integral basis. We call $\mathcal O_{K_\omega}:=\{x\in K_\omega:|x|_\omega\le1\}$ the \emph{valuation ring} of $K_\omega$. The closed unit ball of $\|\mathord{\cdot}\|_{L,\omega}$ is $\mathcal O_{K_\omega}\otimes_\mathbb{Z} L$. Every $\mathbb{Z}$-linear automorphism of $L$ therefore induces an isometry at every non-Archimedean place. The resulting norms are ultrametric or induced by inner products, respectively.

To verify measurability and domination, fix such an integral basis and compare with its reference norm family from \Cref{subsec:norm-families}. The norms agree at every non-Archimedean place. At Archimedean places the Gram matrix of the given metric, and the inverse matrix for the dual metric, are fixed positive-definite real matrices. For each $(x_1,\ldots,x_r)\in K^r$, their quadratic forms can be expressed using the functions $\omega\mapsto|x_i|_\omega^2$ and $\omega\mapsto|x_i+x_j|_\omega^2$ for $(i,j)\in\{1,\ldots,r\}^2$. The eigenvalue bounds give a uniform bound for the distance from the reference norm on $\Omega_\infty$. Since $\nu(\Omega_\infty)$ is finite, \eqref{eq:domination-distance} proves domination. Thus $\overline L^{\,S}$ is a Hermitian adelic vector bundle.

The determinant of the chosen integral basis has norm one at every non-Archimedean place and norm $\operatorname{covol}(L)$ at every Archimedean place. The degree and slope formulas are therefore
\begin{equation}\label{eq:lattice-extension-slope}
	\begin{aligned}
		\operatorname{\widehat{deg}}(\overline L^{\,S}) & =-\nu(\Omega_\infty)\ln(\operatorname{covol}(L)),        \\
		\operatorname{\widehat{\mu}}(\overline L^{\,S})
		                                                & =-\frac{\nu(\Omega_\infty)}r\ln(\operatorname{covol}(L))
		=\nu(\Omega_\infty)\operatorname{\mu}(L).
	\end{aligned}
\end{equation}
On the standard adelic curve of $\mathbb{Q}$, the factor $\nu(\Omega_\infty)$ equals one, so these agree with the degree and slope for Euclidean lattices. We omit the hats for Euclidean lattices and apply \Cref{prop:finitegroup} through this identification.

\section{Rectangular Specht lattices}\label{sec:specht}

The tensor realization below uses the anti-symmetrization and constructions from invariant theory recalled in Appendix A of Chen--Moriwaki (cf. \cite{cmp}). We specify the integral Specht lattice and its metric, and compare the partition and norm conventions with \cite[Sections A.1.2--A.1.4 and A.2.6]{cmp}.

\subsection{Young diagrams and standard tableaux}

Let $(r,n)\in\mathbb{Z}_{\ge1}^2$. We call a tuple $\lambda=(\lambda_1,\ldots,\lambda_r)\in\mathbb{Z}_{\ge1}^r$ a \emph{partition} of $n$ if
\[ \lambda_1\ge\cdots\ge\lambda_r, \qquad\sum_{i=1}^r\lambda_i=n. \]
Define its \emph{Young diagram} by
\[ D_\lambda:=\{(i,j)\in\{1,\ldots,r\}\times\mathbb{Z}_{\ge1}: j\in\{1,\ldots,\lambda_i\}\}. \]
We draw the box $(i,j)$ in row $i$ from the top and column $j$ from the left. For example,
\[ (3,2):\quad\ydiagram{3,2} \qquad\qquad (2,2):\quad\ydiagram{2,2}\;. \]
We call a bijection $t:D_\lambda\to\{1,\ldots,n\}$ a \emph{Young tableau} of shape $\lambda$. Define $\operatorname{YT}(\lambda)$ to be the set of all such tableaux. We call $t\in\operatorname{YT}(\lambda)$ a \emph{standard tableau} if
\begin{gather*}
	\forall\,(i,j)\in D_\lambda\text{ with }(i,j+1)\in D_\lambda,\quad
	t(i,j)<t(i,j+1),\\
	\forall\,(i,j)\in D_\lambda\text{ with }(i+1,j)\in D_\lambda,\quad
	t(i,j)<t(i+1,j).
\end{gather*}
Thus entries increase along rows and down columns. Define $\operatorname{SYT}(\lambda)$ to be the subset of $\operatorname{YT}(\lambda)$ consisting of standard tableaux.

For the partition $\lambda=(\lambda_1,\ldots,\lambda_r)$, define its \emph{conjugate partition} by
\[
	\begin{gathered}
		\widetilde\lambda:=(\widetilde\lambda_1,\ldots,
		\widetilde\lambda_{\lambda_1}),\\
		\forall\,j\in\{1,\ldots,\lambda_1\},\quad
		\widetilde\lambda_j:=\#\{i\in\{1,\ldots,r\}:\lambda_i\ge j\}.
	\end{gathered}
\]
Thus $\widetilde\lambda_j$ is the number of boxes in column $j$ of $D_\lambda$. Its Young diagram is obtained by interchanging rows and columns:
\[ D_{\widetilde\lambda} =\{(j,i):(i,j)\in D_\lambda\}. \]
For $t\in\operatorname{YT}(\lambda)$, define its \emph{transposed tableau} $\widetilde t\in\operatorname{YT}(\widetilde\lambda)$ by
\[ \forall\,(i,j)\in D_\lambda,\quad \widetilde t(j,i):=t(i,j). \]
Transposition interchanges the row and column inequalities, so $t$ is standard if and only if $\widetilde t$ is standard. For example, the partition $(3,2)$ has conjugate partition $(2,2,1)$, and transposition gives
\[ \ytableaushort{123,45} \quad\longmapsto\quad \ytableaushort{14,25,3}\;. \]

For the partition $\lambda=(2,2)$, the two standard tableaux are
\begin{equation}\label{eq:tableaux-square}
	t_1:=\ytableaushort{12,34}\;, \qquad\qquad t_2:=\ytableaushort{13,24}\;.
\end{equation}

\subsection{Tensor words and standard Specht vectors} \label{subsec:specht-tensors}

For the remainder of \Cref{sec:specht}, fix $(a,k)\in\mathbb{Z}_{\ge1}^2$ and define
\[ d:=ak,\qquad \lambda:=(k,\ldots,k)\in\mathbb N^a. \]
Thus $a$ is the number of rows of $D_\lambda$, $k$ is its number of columns, and $d$ is its number of boxes.

\begin{definition}\label{def:tableau-tensor-word}
	Let $(e_1,\ldots,e_a)$ be the standard basis of $\mathbb{Z}^a$. For $b=(b_1,\ldots,b_d)\in\{1,\ldots,a\}^d$, define the \emph{tensor word}
	\[ e_b:=e_{b_1}\otimes\cdots\otimes e_{b_d}. \]
	Equip $(\mathbb{R}^a)^{\otimes d}$ with the inner product for which $(e_b)_{b\in\{1,\ldots,a\}^d}$ is an orthonormal basis. For $t\in\operatorname{YT}(\lambda)$, define $\rho_t:\{1,\ldots,d\}\to\{1,\ldots,a\}$ by
	\[ \forall\,(i,j)\in D_\lambda,\quad \rho_t(t(i,j)):=i, \]
	and define its associated word by
	\begin{equation}\label{eq:tableau-tensor-word}
		w_t:=e_{(\rho_t(1),\ldots,\rho_t(d))}.
	\end{equation}
	Thus the letter in position $j$ records the row containing $j$, and each letter $e_i$ occurs exactly $k$ times.
\end{definition}

\begin{example}\label{ex:tableau-row-reading}
	For $(a,k)=(2,2)$, the tableaux in \eqref{eq:tableaux-square} give the following row maps:
	\[
		\begin{array}{c|cccc}
			j             & 1 & 2 & 3 & 4 \\\hline
			\rho_{t_1}(j) & 1 & 1 & 2 & 2 \\
			\rho_{t_2}(j) & 1 & 2 & 1 & 2
		\end{array}
	\]
	Consequently,
	\begin{align*}
		w_{t_1} & =e_{(1,1,2,2)}=e_1\otimes e_1\otimes e_2\otimes e_2, \\
		w_{t_2} & =e_{(1,2,1,2)}=e_1\otimes e_2\otimes e_1\otimes e_2.
	\end{align*}
	For example, the third factor of $w_{t_1}$ is $e_2$ because the entry $3$ occurs in the second row of $t_1$.
\end{example}

\begin{definition}\label{def:tableau-permutation-actions}
	For $n\in\mathbb{Z}_{\ge1}$, let $\mathfrak S_n$ denote the symmetric group of $\{1,\ldots,n\}$. Define the action of $\sigma\in\mathfrak S_d$ on tensor words by
	\[ \forall\,b\in\{1,\ldots,a\}^d,\quad \sigma\cdot e_b:=e_{(b_{\sigma^{-1}(1)},\ldots,b_{\sigma^{-1}(d)})}, \]
	and extend it linearly. On tableaux, define the relabelling action by
	\[ \forall\,(\sigma,t,(i,j))\in\mathfrak S_d\times\operatorname{YT}(\lambda)\times D_\lambda,\quad (\sigma\cdot t)(i,j):=\sigma(t(i,j)). \]
\end{definition}

\begin{remark}
	The permutation action on tensor words is isometric because it permutes an orthonormal basis. For every $(\sigma,t)\in\mathfrak S_d\times\operatorname{YT}(\lambda)$, the definitions give $w_{\sigma\cdot t}=\sigma\cdot w_t$.
\end{remark}

\begin{example}
	For $(a,k)=(2,2)$, the transposition $(2\;3)$ gives
	\[ (2\;3)\cdot\ytableaushort{12,34}=\ytableaushort{13,24}\;, \qquad (2\;3)\cdot e_{(1,1,2,2)}=e_{(1,2,1,2)}. \]
	Thus $(2\;3)\cdot w_{t_1}=w_{t_2}$. The cycle $(1\;2\;3)$ gives
	\[ (1\;2\;3)\cdot\ytableaushort{12,34}=\ytableaushort{23,14}\;, \qquad (1\;2\;3)\cdot e_{(1,1,2,2)}=e_{(2,1,1,2)}. \]
	Here the factor originally in position $3$ moves to position $1$, which explains the inverse permutation in the definition of the action. Permuting two labels in the same row can leave the word unchanged:
	\[ (1\;2)\cdot\ytableaushort{12,34}=\ytableaushort{21,34}\;, \qquad (1\;2)\cdot e_{(1,1,2,2)}=e_{(1,1,2,2)}. \]

\end{example}

\begin{definition}\label{def:integral-specht-lattice}
	Let $t\in\operatorname{YT}(\lambda)$. For $j\in\{1,\ldots,k\}$, define the column set
	\[ B_j(t):=\{t(i,j):i\in\{1,\ldots,a\}\}, \]
	and define the \emph{column stabilizer} by
	\[ C_t:=\{\sigma\in\mathfrak S_d: \forall\,j\in\{1,\ldots,k\},\ \sigma(B_j(t))=B_j(t)\}. \]
	For $n\in\mathbb{Z}_{\ge1}$, let $\operatorname{sgn}:\mathfrak S_n\to\{-1,1\}$ denote the sign homomorphism. Define the \emph{Specht vector} associated with $t$ by
	\begin{equation}\label{eq:specht-vector-definition}
		v_t:=\sum_{\sigma\in C_t} \operatorname{sgn}(\sigma)\,\sigma\cdot w_t.
	\end{equation}
	For $t\in\operatorname{SYT}(\lambda)$, we call $v_t$ a \emph{standard Specht vector}. Define the \emph{integral Specht lattice} by
	\begin{equation}\label{eq:specht-tensor-lattice}
		L_{a,d}:=\sum_{t\in\operatorname{YT}(\lambda)}\mathbb{Z} v_t \subseteq(\mathbb{Z}^a)^{\otimes d},
	\end{equation}
	with the induced Euclidean metric, and put
	\[ L_{\mathrm{std},a,d}:=\sum_{t\in\operatorname{SYT}(\lambda)}\mathbb{Z} v_t\subseteq L_{a,d}. \]
\end{definition}

\begin{example}\label{ex:one-column-specht-definition}
	Take $(a,k)=(2,1)$, so $d=2$ and $\lambda=(1,1)$. The two tableaux of this shape are
	\[ t:=\ytableaushort{1,2}\;,\qquad u:=\ytableaushort{2,1}\;. \]
	Write $\operatorname{Id}$ for the identity permutation. Their column sets and column stabilizers are
	\[ B_1(t)=B_1(u)=\{1,2\},\qquad C_t=C_u=\{\operatorname{Id},(1\;2)\}. \]
	Their associated tensor words satisfy
	\[ w_t=e_1\otimes e_2,\qquad w_u=e_2\otimes e_1. \]
	Since $\operatorname{sgn}(\operatorname{Id})=1$ and $\operatorname{sgn}((1\;2))=-1$, formula \eqref{eq:specht-vector-definition} gives
	\begin{align*}
		v_t & =w_t-(1\;2)\cdot w_t=e_1\otimes e_2-e_2\otimes e_1,      \\
		v_u & =w_u-(1\;2)\cdot w_u=e_2\otimes e_1-e_1\otimes e_2=-v_t.
	\end{align*}
	Here $\operatorname{YT}(\lambda)=\{t,u\}$ and $\operatorname{SYT}(\lambda)=\{t\}$. Thus the two sums in \Cref{def:integral-specht-lattice} become
	\[ L_{2,2}=\mathbb{Z}v_t+\mathbb{Z}v_u=\mathbb{Z}v_t, \qquad L_{\mathrm{std},2,2}=\mathbb{Z}v_t. \]
	In this example the two lattices coincide. The two tensor words in $v_t$ are orthonormal, so
	\[ \|v_t\|^2=2,\qquad\operatorname{covol}(L_{2,2})=\sqrt2. \]
\end{example}

\begin{remark}\label{rem:standard-specht-normalization}
	The \emph{standard basis theorem}, in the rational form used here, states that
	\[ L_{a,d}\otimes_\mathbb{Z}\mathbb{Q} =\bigoplus_{t\in\operatorname{SYT}(\lambda)}\mathbb{Q} v_t. \]
	Thus every element of $L_{a,d}\otimes_\mathbb{Z}\mathbb{Q}$ has a unique expression as a rational linear combination of the standard Specht vectors. This is the tensor realization of the standard polytabloid basis theorem stated in \cite[Section 2.1, p. 11, paragraph beginning ``A standard basis'']{kn}; see also \cite[Section 2, first two paragraphs]{rh} for its rational form. In particular, $L_{\mathrm{std},a,d}$ has finite index in $L_{a,d}$, and both lattices have rank $\#(\operatorname{SYT}(\lambda))$.

	These are the classical polytabloids in their tensor realization. The map sending a row tabloid to the tensor word whose $j$th letter records the row containing $j$ is a bijection onto the words with $k$ occurrences of each letter. It preserves the orthonormal bases and sends column anti-symmetrization to \eqref{eq:specht-vector-definition}, so the lattice and its metric have the classical Specht normalization. For every $(t,\sigma)\in\operatorname{YT}(\lambda)\times\mathfrak S_d$, $\sigma v_t=v_{\sigma t}$; hence $L_{a,d}$ is stable under $\mathfrak S_d$.
\end{remark}

\begin{definition}\label{def:standard-specht-gram}
	Let $V$ be a Euclidean space, let $m\in\mathbb{Z}_{\ge1}$, and let $(x_1,\ldots,x_m)\in V^m$. Define their \emph{Gram matrix} by
	\[ \operatorname{Gram}(x_1,\ldots,x_m) :=(\langle x_i,x_j\rangle)_{(i,j)\in\{1,\ldots,m\}^2}. \]
	In a fixed ordering of the standard tableaux, put
	\[ G_{\mathrm{std}}:=(\langle v_t,v_u\rangle)_{(t,u)\in\operatorname{SYT}(\lambda)\times\operatorname{SYT}(\lambda)}. \]
\end{definition}

\begin{remark}
	Since the standard Specht vectors form a $\mathbb{Z}$-basis of $L_{\mathrm{std},a,d}$, the determinant of $G_{\mathrm{std}}$ is $\operatorname{covol}(L_{\mathrm{std},a,d})^2$.
\end{remark}

\begin{example}\label{ex:square-specht-lattice}
	Take $(a,k)=(2,2)$, so $d=4$ and $\lambda=(2,2)$. For $t_1$ in \eqref{eq:tableaux-square}, the column sets are $B_1(t_1)=\{1,3\}$ and $B_2(t_1)=\{2,4\}$. Then
	\[ C_{t_1}=\{\operatorname{Id},(1\;3),(2\;4),(1\;3)(2\;4)\}. \]
	The four tableaux in its column orbit are
	\[
		\begin{array}{c@{\qquad}c@{\qquad}c@{\qquad}c}
			t_1                   & (1\;3)\cdot t_1          & (2\;4)\cdot t_1 & (1\;3)(2\;4)\cdot t_1 \\[5pt]
			\ytableaushort{12,34} & \ytableaushort{32,14}    &
			\ytableaushort{14,32} & \ytableaushort{34,12}\;.
		\end{array}
	\]
	Formula \eqref{eq:specht-vector-definition} gives
	\begin{align}
		v_{t_1} & =e_{(1,1,2,2)}-e_{(2,1,1,2)}
		-e_{(1,2,2,1)}+e_{(2,2,1,1)},
		\label{eq:square-specht-one}           \\
		v_{t_2} & =e_{(1,2,1,2)}-e_{(2,1,1,2)}
		-e_{(1,2,2,1)}+e_{(2,1,2,1)}.
		\label{eq:square-specht-two}
	\end{align}
	To compute the full lattice, consider also the nonstandard tableau
	\[ t_3:=\ytableaushort{12,43}\;,\qquad B_1(t_3)=\{1,4\},\quad B_2(t_3)=\{2,3\}. \]
	Its column stabilizer is
	\[ C_{t_3}=\{\operatorname{Id},(1\;4),(2\;3),(1\;4)(2\;3)\}. \]
	Applying \eqref{eq:specht-vector-definition} once more gives
	\begin{align*}
		v_{t_3}
		 & =e_{(1,1,2,2)}-e_{(2,1,2,1)}-e_{(1,2,1,2)}+e_{(2,2,1,1)} \\
		 & =v_{t_1}-v_{t_2}.
	\end{align*}
	The column sets of any tableau of shape $(2,2)$ give one of the three unordered pairings
	\[ \{\{1,3\},\{2,4\}\},\qquad \{\{1,2\},\{3,4\}\},\qquad \{\{1,4\},\{2,3\}\}. \]
	These are the pairings for $t_1$, $t_2$, and $t_3$, respectively. Interchanging the two columns leaves both the tensor word and the column stabilizer unchanged, hence leaves the Specht vector unchanged. Interchanging the two entries of one column multiplies the Specht vector by $-1$, by reindexing the alternating sum in \eqref{eq:specht-vector-definition}. Each pairing admits two orders of its columns and two orders of the entries in each column, giving $2\cdot2^2=8$ tableaux. Thus all $24$ tableaux have Specht vector in $\{\pm v_{t_1},\pm v_{t_2},\pm v_{t_3}\}$. It follows that
	\[ L_{2,4}=\mathbb{Z}v_{t_1}+\mathbb{Z}v_{t_2}+\mathbb{Z}v_{t_3} =\mathbb{Z}v_{t_1}+\mathbb{Z}v_{t_2} =L_{\mathrm{std},2,4}. \]
	The coefficients of $e_{(1,1,2,2)}$ and $e_{(1,2,1,2)}$ in \eqref{eq:square-specht-one}--\eqref{eq:square-specht-two} show that $v_{t_1}$ and $v_{t_2}$ are linearly independent. Consequently they form a $\mathbb{Z}$-basis of the full lattice.

	The tensor words are orthonormal. Each of $v_{t_1}$ and $v_{t_2}$ has four coefficients equal to $1$ or $-1$, and their common words are $e_{(2,1,1,2)}$ and $e_{(1,2,2,1)}$, with coefficient $-1$ in both vectors. Therefore
	\[
		\|v_{t_1}\|^2=\|v_{t_2}\|^2=4,\qquad \langle v_{t_1},v_{t_2}\rangle=2,\qquad G_{\mathrm{std}}=\begin{pmatrix}4&2\\2&4\end{pmatrix}.
	\]
	In particular,
	\[ \operatorname{det}(G_{\mathrm{std}})=4\cdot4-2^2=12,\qquad \operatorname{covol}(L_{2,4}) =\operatorname{covol}(L_{\mathrm{std},2,4})=\sqrt{12}=2\sqrt3. \]

\end{example}

\begin{example}
	For $(a,k)=(1,3)$, the partition is $\lambda=(3)$ and its unique standard tableau is
	\[ t_{\mathrm{row}}:=\ytableaushort{123}\;. \]
	Every column contains one entry, so $C_{t_{\mathrm{row}}}=\{\operatorname{Id}\}$ and
	\[ v_{t_{\mathrm{row}}}=w_{t_{\mathrm{row}}} =e_{(1,1,1)}=e_1^{\otimes3},\qquad L_{1,3}=\mathbb{Z} e_1^{\otimes3},\qquad \|v_{t_{\mathrm{row}}}\|=1. \]
	For $(a,k)=(3,1)$, the partition is $\lambda=(1,1,1)$ and its unique standard tableau is
	\[ t_{\mathrm{vert}}:=\ytableaushort{1,2,3}\;. \]
	Here $w_{t_{\mathrm{vert}}}=e_{(1,2,3)}$ and $C_{t_{\mathrm{vert}}}=\mathfrak S_3$. Thus
	\begin{align*}
		v_{t_{\mathrm{vert}}}
		 & =e_{(1,2,3)}-e_{(1,3,2)}-e_{(2,1,3)}       \\
		 & \quad+e_{(2,3,1)}+e_{(3,1,2)}-e_{(3,2,1)}.
	\end{align*}
	The six words are distinct and orthonormal. Consequently,
	\[ L_{3,3}=\mathbb{Z} v_{t_{\mathrm{vert}}},\qquad \|v_{t_{\mathrm{vert}}}\|^2=6,\qquad \operatorname{covol}(L_{3,3})=\sqrt6. \]

\end{example}

\begin{example}
	For $(a,k)=(2,3)$, consider the standard tableau
	\[ t_{\mathrm{alt}}:=\ytableaushort{135,246} \qquad\text{of shape }(3,3). \]
	Its word is $w_{t_{\mathrm{alt}}}=e_{(1,2,1,2,1,2)}$. The column sets are $\{1,2\}$, $\{3,4\}$, and $\{5,6\}$, so the column stabilizer is generated by the three independent transpositions $(1\;2)$, $(3\;4)$, and $(5\;6)$. Antisymmetrizing these three pairs gives
	\[ v_{t_{\mathrm{alt}}} =(e_1\otimes e_2-e_2\otimes e_1)^{\otimes3}. \]
	Expanding the right-hand side gives eight distinct orthonormal words, each with coefficient $1$ or $-1$. Hence
	\[ \|v_{t_{\mathrm{alt}}}\|^2=8. \]

\end{example}

\subsection{Determinant tensors and anti-symmetrization}

Let $R$ be a commutative ring with identity and let $V$ be a free $R$-module of rank $a$. Define $\operatorname{det}(V):=\bigwedge\nolimits^a V$. Let $\operatorname{GL}(V)$ denote the group of linear automorphisms of $V$, and define its \emph{special linear group} by
\[ \operatorname{SL}(V):=\{g\in\operatorname{GL}(V):\operatorname{det}(g)=1\}. \]
Following \cite[Section A.1.2]{cmp}, define the \emph{anti-symmetrization map}, for $(v_1,\ldots,v_a)\in V^a$, by
\begin{equation}\label{eq:determinant-antisymmetrization}
	\iota_a:\operatorname{det}(V)\longrightarrow V^{\otimes a},\qquad \iota_a(v_1\wedge\cdots\wedge v_a):= \sum_{\sigma\in\mathfrak S_a}\operatorname{sgn}(\sigma) v_{\sigma(1)}\otimes\cdots\otimes v_{\sigma(a)}.
\end{equation}
This map is injective and equivariant for the general linear group of $V$, as recalled in that section. For $R=\mathbb{Z}$ and $V=\mathbb{Z}^a$, define the \emph{determinant tensor}
\[ \Delta_a:=\iota_a(e_1\wedge\cdots\wedge e_a) \]
and the \emph{column-reading tableau} by
\[ \forall\,(i,j)\in\{1,\ldots,a\}\times\{1,\ldots,k\},\quad t_{\mathrm{col}}(i,j):=i+a(j-1). \]
Applying \eqref{eq:determinant-antisymmetrization} to each column gives
\begin{equation}\label{eq:column-tableau-determinant}
	v_{t_{\mathrm{col}}}=\Delta_a^{\otimes k}.
\end{equation}
For $(a,k)=(2,2)$, this identity reads
\[ t_{\mathrm{col}}=t_2=\ytableaushort{13,24}\;,\qquad v_{t_2}=(e_1\otimes e_2-e_2\otimes e_1)^{\otimes2}. \]
Every $t\in\operatorname{YT}(\lambda)$ is obtained from $t_{\mathrm{col}}$ by relabelling. {The definition of the full lattice therefore gives}
\begin{equation}\label{eq:specht-determinant-orbit}
	L_{a,d}=\sum_{\sigma\in\mathfrak S_d} \mathbb{Z}\,\sigma\cdot\Delta_a^{\otimes k}.
\end{equation}
The group $\operatorname{SL}(\mathbb{Z}^a)$ fixes $\Delta_a$ and commutes with permutations of tensor positions. Hence every element of $L_{a,d}$ is fixed by $\operatorname{SL}(\mathbb{Z}^a)$.

\subsection{Invariant theory and Schur functor conventions}

{To identify the full invariant space, assume that $\operatorname{char}(K)=0$, and put $V:=K^a$.} For $m\in\mathbb{Z}_{\ge1}$, let $S^m(V)$ denote the \emph{symmetric power}, defined as the quotient of $V^{\otimes m}$ by the $K$-submodule generated by the elements $x-\sigma x$ for $(x,\sigma)\in V^{\otimes m}\times\mathfrak S_m$, as in \cite[Section A.1.1]{cmp}. Apply the first fundamental theorem (cf. \cite[Theorem A.1.1]{cmp}) to $d$ labelled copies of $V$ and take the component of degree one in each copy. This component identifies with $V^{\otimes d}$, and its invariant vectors are spanned by products of $k$ determinant tensors on disjoint sets of $a$ positions. These products are precisely the permutations of $\Delta_a^{\otimes k}$. Hence
\begin{equation}\label{eq:specht-invariant-space}
	{L_{a,d}\otimes_\mathbb{Z} K} =(V^{\otimes d})^{\operatorname{SL}(V)}.
\end{equation}

\begin{proposition}\label{prop:specht-absolute-irreducibility}
	The permutation representation of $\mathfrak S_d$ on $L_{a,d}\otimes_\mathbb{Z}K$ is absolutely irreducible.
\end{proposition}
\begin{proof}
	Put $U:=L_{a,d}\otimes_\mathbb{Z}\mathbb{Q}$, and let $\rho:\mathfrak S_d\to\operatorname{GL}(U)$ denote the permutation representation. For a field extension $K'/\mathbb{Q}$, define $U_{K'}:=K'\otimes_\mathbb{Q}U\simeq L_{a,d}\otimes_\mathbb{Z}K'$ and denote the extended representation by $\rho_{K'}$. We will prove irreducibility of $U_{K'}$ for every such extension. Taking $K'$ to be any extension of $K$ then proves the assertion.

	First consider $K'=\mathbb{C}$. We specify the connection with \cite[formula (3.2) and the proof of Lemma 3.9, p. 17]{bgow}, where the representations are written as spaces of linear forms. The complex-linear map sending each tensor word in $(\mathbb{C}^a)^{\otimes d}$ to its coordinate functional is an isomorphism with the dual tensor space. Since permutations permute the basis of tensor words, this map intertwines the permutation action with the contragredient action on the dual. It sends $\Delta_a^{\otimes k}$ to the product of determinant forms in \cite[formula (3.2)]{bgow}, with their parameters $n=a$, $m=d$, and $\pi=\operatorname{Id}$. Thus \eqref{eq:specht-determinant-orbit} identifies $U_{\mathbb{C}}$ with the span of the permutation orbit of that product of determinant forms. The proof of \cite[Lemma 3.9, p. 17]{bgow} shows, by Schur--Weyl duality, that this span is an irreducible $\mathfrak S_d$-representation. Hence $U_{\mathbb{C}}$ is irreducible.

	Burnside's irreducibility criterion (cf. \cite[Theorem 2.7.24, p. 98]{kowalski}) states that, for a finite-dimensional irreducible representation over an algebraically closed field, the linear span of the representing operators is the full endomorphism algebra. Applied over $\mathbb{C}$, it gives
	\begin{equation}\label{eq:specht-complex-action-span}
		\sum_{\sigma\in\mathfrak S_d}\mathbb{C}\rho_{\mathbb{C}}(\sigma) =\operatorname{End}_{\mathbb{C}}(U_{\mathbb{C}}).
	\end{equation}
	To descend this equality to $\mathbb{Q}$, define
	\[ A:=\sum_{\sigma\in\mathfrak S_d}\mathbb{Q}\rho(\sigma) \subseteq\operatorname{End}_{\mathbb{Q}}(U). \]
	This span is a subalgebra, since $\rho(\sigma)\rho(\tau)=\rho(\sigma\tau)$ for every $(\sigma,\tau)\in\mathfrak S_d\times\mathfrak S_d$. For every field extension $K'/\mathbb{Q}$, the natural map
	\[ K'\otimes_\mathbb{Q}\operatorname{End}_{\mathbb{Q}}(U) \longrightarrow\operatorname{End}_{K'}(U_{K'}) \]
	is an isomorphism. Indeed, after choosing any $\mathbb{Q}$-basis of $U$ and the induced $K'$-basis of $U_{K'}$, it extends each rational matrix to the same matrix over $K'$. In particular, it sends the basis of matrix units to the basis of matrix units. Extending a $\mathbb{Q}$-basis of $A$ to a basis of $\operatorname{End}_{\mathbb{Q}}(U)$ also shows that the induced map on $K'\otimes_\mathbb{Q}A$ is injective. Its image is exactly the $K'$-span of the operators $\rho_{K'}(\sigma)$, because $A$ is generated as a vector space by the operators $\rho(\sigma)$. Consequently \eqref{eq:specht-complex-action-span} gives
	\[
		\begin{aligned}
			\operatorname{dim}_{\mathbb{Q}}(A)
			 & =\operatorname{dim}_{\mathbb{C}}(\mathbb{C}\otimes_\mathbb{Q}A) \\
			 & =\operatorname{dim}_{\mathbb{C}}
			(\operatorname{End}_{\mathbb{C}}(U_{\mathbb{C}}))                  \\
			 & =(\operatorname{dim}_{\mathbb{Q}}(U))^2
			=\operatorname{dim}_{\mathbb{Q}}
			(\operatorname{End}_{\mathbb{Q}}(U)).
		\end{aligned}
	\]
	Since $A$ is a subspace of $\operatorname{End}_{\mathbb{Q}}(U)$, equality of these finite dimensions implies $A=\operatorname{End}_{\mathbb{Q}}(U)$. Extending scalars once more therefore yields, for every field extension $K'/\mathbb{Q}$,
	\begin{equation}\label{eq:specht-action-span}
		\sum_{\sigma\in\mathfrak S_d}K'\rho_{K'}(\sigma) =\operatorname{End}_{K'}(U_{K'}).
	\end{equation}

	Finally, let $W\subseteq U_{K'}$ be a nonzero $\mathfrak S_d$-invariant subspace. It is preserved by every linear combination of the operators $\rho_{K'}(\sigma)$, hence by every endomorphism of $U_{K'}$, by \eqref{eq:specht-action-span}. Choose $u\in W\setminus\{0\}$. For every $v\in U_{K'}$, extending $u$ to a basis gives an endomorphism $T$ with $T(u)=v$. Thus $v=T(u)\in W$, so $W=U_{K'}$. This proves the required irreducibility over every extension of $\mathbb{Q}$.
\end{proof}

For the rectangular partition $\lambda=(k,\ldots,k)\in\mathbb N^a$, each of its $k$ columns has $a$ boxes. Hence its conjugate partition is
\[ \widetilde\lambda=(a,\ldots,a)\in\mathbb N^k. \]
We use this conjugate partition to compare the Schur functor conventions. Partitions are extended by zero coordinates when using the sequence notation of \cite[Section A.1.3]{cmp}. Let $r\in\mathbb{Z}_{\ge1}$ and let $\mu=(\mu_1,\ldots,\mu_r)$ be a partition of $d$. Following that section, define the \emph{Schur module} $L^\mu(V)$ to be the image of the composite
\[ \bigotimes_{i=1}^r\bigwedge\nolimits^{\mu_i}V \longrightarrow V^{\otimes d} \longrightarrow \bigotimes_{j=1}^{\mu_1}S^{\widetilde\mu_j}(V). \]
The first arrow applies unnormalized anti-symmetrization to each exterior power and orders the tensor factors along the rows of $D_\mu$. The second arrow groups these factors by columns and then takes their symmetric products. Explicitly, for $(v_{ij})_{(i,j)\in D_\mu}\in V^{D_\mu}$, the second arrow sends
\[ \bigotimes_{i=1}^r\bigotimes_{j=1}^{\mu_i}v_{ij} \longmapsto \bigotimes_{j=1}^{\mu_1} (v_{1j}\cdots v_{\widetilde\mu_j,j}), \]
where the product in the $j$-th factor denotes the class in $S^{\widetilde\mu_j}(V)$. Linear maps between finite-dimensional $K$-vector spaces induce maps between the corresponding images. We call the resulting functor $V\mapsto L^\mu(V)$ the \emph{Schur functor} in this convention (cf. \cite[Section A.1.3, p. 221]{cmp}). For example, $L^{(d)}(V)=\bigwedge\nolimits^dV$, whereas $\mu=(1,\ldots,1)\in\mathbb N^d$ gives $L^\mu(V)=S^d(V)$.

For $\mu=\widetilde\lambda$, the defining composite becomes
\[ \operatorname{det}(V)^{\otimes k} \xrightarrow{\ \iota_a^{\otimes k}\ }V^{\otimes d} \longrightarrow S^k(V)^{\otimes a}. \]
Here the columns of $D_\mu$ are the rows of $D_\lambda$. The second arrow therefore groups positions by rows of $D_\lambda$ and passes to the symmetric power in each row. For $i\in\{1,\ldots,a\}$, write $e_i^k$ for the class of $e_i^{\otimes k}$ in $S^k(V)$. The image of the displayed map, $L^{\widetilde\lambda}(V)$, is isomorphic to the one-dimensional representation $\operatorname{det}(V)^{\otimes k}$: the map is equivariant, and the coefficient of $e_1^k\otimes\cdots\otimes e_a^k$ in the image of $(e_1\wedge\cdots\wedge e_a)^{\otimes k}$ is one. The dimension of ${L_{a,d}\otimes_\mathbb{Z} K}$ is the multiplicity of this representation in $V^{\otimes d}$, and the commuting $\mathfrak S_d$-action acts on its multiplicity space. Thus the Schur functor indexed by $\widetilde\lambda$ and the Specht lattice indexed by $\lambda$ play different roles.

For a finite family of finite-dimensional vector spaces over a perfect field, \cite[Theorem A.1.4, formula (A.1.4)]{cmp} states that the vectors in a symmetric power of their tensor product invariant under the product of their special linear groups are spanned by the images of products of determinant tensors, after independent permutations of tensor positions in the factors. This is the lifting statement for several factors associated with the construction above. It supplies an algebraic spanning statement; the metric normalization is specified separately below.

Let $n\in\mathbb{Z}_{\ge1}$, let $\nu$ be a partition of $n$, and let $(i,j)\in D_\nu$. Define the \emph{hook} at $(i,j)$ by
\[ H_{ij}:=\{(i,q)\in D_\nu:q\ge j\} \cup\{(p,j)\in D_\nu:p>i\}, \]
and define its \emph{hook length} by
\[ h_{ij}:=\#(H_{ij})=\nu_i-j+\widetilde\nu_j-i+1. \]
The hook consists of the box $(i,j)$, the boxes to its right in the same row, and the boxes below it in the same column. For example, for $\nu=(3,2)$, the following diagram records the hook length in each box:
\[ \ytableaushort{431,21}\;. \]
Here $H_{11}=\{(1,1),(1,2),(1,3),(2,1)\}$, so $h_{11}=4$. The \emph{hook-length formula} of Frame--Robinson--Thrall (cf. \cite{frt}) states that
\[ \#(\operatorname{SYT}(\nu))=\frac{n!}{\prod_{(i,j)\in D_\nu}h_{ij}}. \]
In this example it gives $\#(\operatorname{SYT}((3,2)))=5!/(4\cdot3\cdot1\cdot2\cdot1)=5$.

\begin{proposition}\label{prop:specht-rank}
	Let $(a,k)\in\mathbb{Z}_{\ge1}^2$, put $d:=ak$, and define $\lambda:=(k,\ldots,k)\in\mathbb N^a$. The integral Specht lattice of \eqref{eq:specht-tensor-lattice} has rank
	\begin{equation}\label{eq:dimension}
		\operatorname{rk}(L_{a,d})=\#(\operatorname{SYT}(\lambda)) =\frac{d!\prod_{j=1}^{a-1}j!} {\prod_{j=0}^{a-1}(k+j)!}.
	\end{equation}
	Empty products in this formula are taken to be one.
\end{proposition}
\begin{proof}
	The rational form of the standard basis theorem stated in \Cref{subsec:specht-tensors}, from \cite[Section 2.1, p. 11]{kn}, gives
	\[ \operatorname{rk}(L_{a,d})=\operatorname{dim}_\mathbb{Q}(L_{a,d}\otimes_\mathbb{Z}\mathbb{Q})=\#(\operatorname{SYT}(\lambda)). \]
	To evaluate this number, let $n\in\mathbb{Z}_{\ge1}$ and let $\nu$ be a partition of $n$ with at most $a$ rows. Append zero parts to obtain $(\nu_1,\ldots,\nu_a)$, and define
	\[ \forall\,i\in\{1,\ldots,a\},\quad b_i:=\nu_i+a-i. \]
	Then $b_1>\cdots>b_a\ge0$. We claim that, for each $i\in\{1,\ldots,a\}$,
	\[ \{h_{ij}:j\in\{1,\ldots,\nu_i\}\} =\{1,\ldots,b_i\} \setminus\{b_i-b_r:r\in\{i+1,\ldots,a\}\}. \]
	Indeed, the hook lengths in row $i$ are positive, are at most $b_i$, and are pairwise distinct, since
	\[ \forall\,j\in\{1,\ldots,\nu_i-1\},\quad h_{ij}-h_{i,j+1}=1+\widetilde\nu_j-\widetilde\nu_{j+1}>0. \]
	For $r\in\{i+1,\ldots,a\}$, the integers $b_i-b_r$ are also distinct and belong to $\{1,\ldots,b_i\}$. For $(j,r)\in\{1,\ldots,\nu_i\}\times\{i+1,\ldots,a\}$, one has
	\[ h_{ij}-(b_i-b_r)=\nu_r-j+\widetilde\nu_j-r+1. \]
	If $r\le\widetilde\nu_j$, then $\nu_r\ge j$, so the right-hand side is at least one. If $r>\widetilde\nu_j$, then $\nu_r<j$, so it is at most minus one. Hence the two sets of integers are disjoint. Their cardinalities add up to $\nu_i+(a-i)=b_i$, proving the claim, including the case $\nu_i=0$, when the set of hooks is empty. Taking products gives
	\[ \prod_{j=1}^{\nu_i}h_{ij} =\frac{b_i!}{\prod_{r=i+1}^a(b_i-b_r)}. \]
	Multiplication over $i\in\{1,\ldots,a\}$, followed by the hook-length formula, yields
	\begin{equation}\label{eq:shifted-count}
		\#(\operatorname{SYT}(\nu))=\frac{n!\prod_{1\le i<j\le a}(b_i-b_j)}{\prod_{i=1}^a b_i!}.
	\end{equation}
	For $\nu=\lambda=(k,\ldots,k)$, one has $n=d$, $b_i=k+a-i$, and $b_i-b_j=j-i$. Therefore
	\begin{align*}
		\prod_{1\le i<j\le a}(b_i-b_j)
		 & =\prod_{j=2}^a\prod_{i=1}^{j-1}(j-i)
		=\prod_{j=2}^a(j-1)!
		=\prod_{j=1}^{a-1}j!,                   \\
		\prod_{i=1}^a b_i!
		 & =\prod_{i=1}^a(k+a-i)!
		=\prod_{r=0}^{a-1}(k+r)!.
	\end{align*}
	Substitution into \eqref{eq:shifted-count} proves \eqref{eq:dimension}. The same counting formulas include the empty partition, with the conventions recalled in \Cref{sec:seminormal}.
\end{proof}

\subsection{Norm normalization and covolume}

Let $\omega\in\Omega$ and let $W$ be a normed $K_\omega$-vector space of dimension $a$. Assume that its norm is ultrametric if $\omega\notin\Omega_\infty$, and induced by an inner product if $\omega\in\Omega_\infty$. Equip $\operatorname{det}(W)$ with the determinant norm and $W^{\otimes a}$ with the corresponding $\varepsilon$-tensor or Hilbert--Schmidt tensor norm. The precise norm comparison for \eqref{eq:determinant-antisymmetrization} is
\begin{equation}\label{eq:antisymmetrization-norm}
	\forall\,\eta\in\operatorname{det}(W),\quad \|\iota_a(\eta)\|:= \begin{cases}
		\|\eta\|_{\omega,\operatorname{det}},            & \text{if the absolute value is non-Archimedean}, \\
		\sqrt{a!}\,\|\eta\|_{\omega,\operatorname{det}}, & \text{if it is Archimedean},
	\end{cases}
\end{equation}
The first equality is \cite[Proposition A.2.8]{cmp}; the second is \cite[Remark A.2.9]{cmp}. In our Euclidean normalization, \eqref{eq:column-tableau-determinant} therefore gives
\begin{equation}\label{eq:column-tableau-norm}
	\|v_{t_{\mathrm{col}}}\|^2 =\|\Delta_a^{\otimes k}\|^2=(a!)^k.
\end{equation}
{We retain this factor in the Specht metric. \Cref{sec:seminormal,sec:asymptotic} will estimate the covolume using the whole Specht basis.}

By \eqref{eq:lattice-extension-slope}, the notation of \Cref{sec:introduction} gives
\[ -\operatorname{\widehat{\mu}}(\overline L_{a,d}^{\,S}) =\nu(\Omega_\infty)\frac{\ln(\operatorname{covol}(L_{a,d}))}{\#(\operatorname{SYT}(\lambda))}. \]
The inclusion $L_{\mathrm{std},a,d}\subseteq L_{a,d}$ gives $\operatorname{covol}(L_{a,d})\le\operatorname{covol}(L_{\mathrm{std},a,d})$. For any $\mathbb{Z}$-basis of $L_{a,d}$, the Gram matrix is positive definite with integer entries, so its determinant is a positive integer. Recall from \Cref{subsec:specht-tensors} that $G_{\mathrm{std}}$ denotes the Gram matrix of the standard Specht vectors $(v_t)_{t\in\operatorname{SYT}(\lambda)}$. Thus
\begin{equation}\label{eq:covol-gram-bound}
	0\le\frac{2\ln(\operatorname{covol}(L_{a,d}))}{\#(\operatorname{SYT}(\lambda))} \le\frac{\ln\operatorname{det} G_{\mathrm{std}}}{\#(\operatorname{SYT}(\lambda))}.
\end{equation}

By \Cref{prop:finitegroup}, $L_{a,d}$ is semistable. Any tensor product of these lattices, and its dual, is semistable under the independent actions of the permutation groups. The diagonal action used in later contractions is distinct from these independent actions.

\section{An exact seminormal norm formula}\label{sec:seminormal}

The local transposition formula used in this section is the classical formula of \cite[Section 2, equation (1)]{rh}. We express its norm ratios in terms of successive row lengths and verify the initial value in our tensor metric. This verification is needed because rescaling every basis vector changes the covolume by that scale raised to the rank.

\begin{definition}\label{def:partition-conventions}
	We define the \emph{empty partition} of zero to be the empty tuple, denoted by $\varnothing$. We identify finite tuples differing only by trailing zeros. For a partition $\mu=(\mu_1,\ldots,\mu_\ell)$ with $\ell\in\mathbb{Z}_{\ge0}$ positive parts, define
	\[ \forall\,i\in\mathbb{Z}_{\ge\ell+1},\quad \mu_i:=0. \]
	This convention applies to all partitions in the remainder of the paper, including prefix shapes and complement partitions. Define $\operatorname{YT}(\varnothing):=\operatorname{SYT}(\varnothing):=\{\varnothing\}$, where the element $\varnothing$ is the unique map from the empty diagram to the empty set.
\end{definition}

\begin{definition}\label{def:partition-containment}
	Let $\lambda$ and $\nu$ be partitions, possibly of different integers. We say that $\nu$ is \emph{contained} in $\lambda$, and write $\nu\subseteq\lambda$, if $ D_\nu\subseteq D_\lambda$, or equivalently,
	\[ \forall\,i\in\mathbb{Z}_{\ge1},\quad \nu_i\le\lambda_i. \]
	Here the coordinates are extended by zero, and equality of the two partitions is allowed.
\end{definition}

\begin{example}
	The partitions $(2,1)$ and $(3,3)$ satisfy $(2,1)\subseteq(3,3)$.
\end{example}

\begin{definition}\label{def:partition-dominance}
	Let $n\in\mathbb{Z}_{\ge0}$ and let $\lambda$ and $\nu$ be partitions of $n$, with their zero extensions as above. We say that $\lambda$ \emph{dominates} $\nu$, and write $\lambda\unrhd\nu$, if
	\[ \forall\,m\in\mathbb{Z}_{\ge1},\quad \sum_{i=1}^m\lambda_i\ge\sum_{i=1}^m\nu_i. \]
\end{definition}
\begin{example}
	The partitions $(3)$ and $(2,1)$ satisfy $(3)\unrhd(2,1)$: their zero-extended tuples begin with $(3,0)$ and $(2,1)$, so the first partial sums satisfy $3\ge2$ and every subsequent partial sum is $3$ on both sides.
\end{example}

\begin{definition}\label{def:tableau-prefix-shape}
	Let $\mu=(\mu_1,\ldots,\mu_\ell)$ be a partition of $n\in\mathbb{Z}_{\ge1}$ and let $t\in\operatorname{SYT}(\mu)$. For $q\in\{0,\ldots,n\}$, define its \emph{prefix shape} by
	\[ \operatorname{sh}_q(t):= \bigl(\#\{j\in\{1,\ldots,\mu_i\}:t(i,j)\le q\}\bigr)_{i=1}^{\ell}. \]
	The row lengths are weakly decreasing, so this tuple represents a partition of $q$. In particular, $\operatorname{sh}_0(t)=\varnothing$ and $\operatorname{sh}_n(t)=\mu$ with the preceding convention.
\end{definition}

To obtain the prefix, retain exactly the boxes whose labels are at most $q$. Since entries increase along rows and down columns, every box above or to the left of a retained box is also retained. Hence these boxes form a Young diagram, and $\operatorname{sh}_q(t)$ records its row lengths.

\begin{example}\label{ex:tableau-prefixes}
	Consider
	\[ t:=\ytableaushort{125,36,4} \in\operatorname{SYT}((3,2,1)). \]
	Its prefixes and their shapes are as follows; the empty prefix is represented by $\varnothing$:
	\[
		\renewcommand{\arraystretch}{1.3} \begin{array}{c|ccccccc}
		q&0&1&2&3&4&5&6\\[3pt]\hline
		\text{prefix}&\varnothing&\hspace{3pt}\raisebox{0pt}[\dimexpr\height+4pt\relax][\dimexpr\depth+4pt\relax]{$\ytableaushort{1}$}\hspace{3pt}&\hspace{3pt}\raisebox{0pt}[\dimexpr\height+4pt\relax][\dimexpr\depth+4pt\relax]{$\ytableaushort{12}$}\hspace{3pt}
		&\hspace{3pt}\raisebox{0pt}[\dimexpr\height+4pt\relax][\dimexpr\depth+4pt\relax]{$\ytableaushort{12,3}$}\hspace{3pt}&\hspace{3pt}\raisebox{0pt}[\dimexpr\height+4pt\relax][\dimexpr\depth+4pt\relax]{$\ytableaushort{12,3,4}$}\hspace{3pt}
		&\hspace{3pt}\raisebox{0pt}[\dimexpr\height+4pt\relax][\dimexpr\depth+4pt\relax]{$\ytableaushort{125,3,4}$}\hspace{3pt}&\hspace{3pt}\raisebox{0pt}[\dimexpr\height+4pt\relax][\dimexpr\depth+4pt\relax]{$\ytableaushort{125,36,4}$}\hspace{3pt}\\[7pt]
		\operatorname{sh}_q(t)&\varnothing&(1)&(2)&(2,1)&(2,1,1)
		&(3,1,1)&(3,2,1)
		\end{array}
	\]
	The shape of the whole tableau does not determine its prefix shapes. For the two tableaux of \eqref{eq:tableaux-square}, the prefixes at $q=2$ give
	\[
		\renewcommand{\arraystretch}{1.3} \begin{array}{c|cc}
		&\hspace{3pt}\raisebox{0pt}[\dimexpr\height+4pt\relax][\dimexpr\depth+4pt\relax]{$t_1=\ytableaushort{12,34}$}\hspace{3pt}&\hspace{3pt}\raisebox{0pt}[\dimexpr\height+4pt\relax][\dimexpr\depth+4pt\relax]{$t_2=\ytableaushort{13,24}$}\hspace{3pt}\\[7pt]\hline
		\text{prefix at }q=2&\hspace{3pt}\raisebox{0pt}[\dimexpr\height+4pt\relax][\dimexpr\depth+4pt\relax]{$\ytableaushort{12}$}\hspace{3pt}&\hspace{3pt}\raisebox{0pt}[\dimexpr\height+4pt\relax][\dimexpr\depth+4pt\relax]{$\ytableaushort{1,2}$}\hspace{3pt}\\[7pt]
		\text{prefix shape}&(2)&(1,1)
		\end{array}
	\]
	Thus both tableaux have shape $(2,2)$, but their two-box prefixes have different shapes.
\end{example}

\begin{definition}\label{def:tableau-dominance}
	Let $n\in\mathbb{Z}_{\ge1}$, and let $t$ and $u$ be standard tableaux with $n$ boxes. We say that $t$ \emph{dominates} $u$, and write $t\unrhd u$, if
	\[ \forall\,q\in\{1,\ldots,n\},\quad \operatorname{sh}_q(t)\unrhd\operatorname{sh}_q(u). \]
	The two prefix shapes at a given $q$ may have different lengths; their comparison always uses the zero extensions. The comparison must hold at every prefix. At a given $q$, the partition dominance condition compares the number of labels at most $q$ in the first $m$ rows, for each $m\in\mathbb{Z}_{\ge1}$.
\end{definition}

\begin{example}\label{ex:tableau-dominance}
	Consider again
	\[ t_1=\ytableaushort{12,34}\;,\qquad t_2=\ytableaushort{13,24}\;. \]
	All their nonempty prefix shapes are listed below:
	\[
		\renewcommand{\arraystretch}{1.3} \begin{array}{c|cccc}
			q                        & 1   & 2     & 3     & 4     \\[3pt]\hline
			\operatorname{sh}_q(t_1) & (1) & (2)   & (2,1) & (2,2) \\
			\operatorname{sh}_q(t_2) & (1) & (1,1) & (2,1) & (2,2)
		\end{array}
	\]
	The shapes agree for $q\in\{1,3,4\}$. At $q=2$, the first row of the first prefix contains two boxes, and that of the second contains one; the sums over the first $m$ rows are both two for every $m\in\mathbb{Z}_{\ge2}$. Hence
	\[ \operatorname{sh}_2(t_1)=(2)\unrhd(1,1) =\operatorname{sh}_2(t_2), \qquad t_1\unrhd t_2. \]
	The reverse dominance fails at $q=2$.
\end{example}

\begin{example}\label{ex:incomparable-tableaux}
	Two standard tableaux of the same shape need not be comparable. Take
	\[ u:=\ytableaushort{125,346}\;,\qquad v:=\ytableaushort{134,256}\;,\qquad (u,v)\in\operatorname{SYT}((3,3))^2. \]
	The following two prefixes already show the obstruction:
	\[
		\renewcommand{\arraystretch}{1.3} \begin{array}{c|cc|cc}
		q&\text{prefix of }u&\text{prefix of }v
		&\operatorname{sh}_q(u)&\operatorname{sh}_q(v)\\[3pt]\hline
		2&\hspace{3pt}\raisebox{0pt}[\dimexpr\height+4pt\relax][\dimexpr\depth+4pt\relax]{$\ytableaushort{12}$}\hspace{3pt}
		&\hspace{3pt}\raisebox{0pt}[\dimexpr\height+4pt\relax][\dimexpr\depth+4pt\relax]{$\ytableaushort{1,2}$}\hspace{3pt}&(2)&(1,1)\\[7pt]
		4&\hspace{3pt}\raisebox{0pt}[\dimexpr\height+4pt\relax][\dimexpr\depth+4pt\relax]{$\ytableaushort{12,34}$}\hspace{3pt}
		&\hspace{3pt}\raisebox{0pt}[\dimexpr\height+4pt\relax][\dimexpr\depth+4pt\relax]{$\ytableaushort{134,2}$}\hspace{3pt}&(2,2)&(3,1)
		\end{array}
	\]
	At $q=2$, the first-row counts give $1<2$, so $\operatorname{sh}_2(v)\ntrianglerighteq\operatorname{sh}_2(u)$. At $q=4$, they give $2<3$, so $\operatorname{sh}_4(u)\ntrianglerighteq\operatorname{sh}_4(v)$. Consequently,
	\[ u\ntrianglerighteq v,\qquad v\ntrianglerighteq u. \]
	This illustrates why one must choose a total order compatible with dominance before applying Gram--Schmidt below.
\end{example}

\begin{definition}\label{def:tableau-content}
	Let $n\in\mathbb{Z}_{\ge1}$ and let $\mu$ be a partition of $n$. For $(t,h)\in\operatorname{SYT}(\mu)\times\{1,\ldots,n\}$, let $(i,j)\in D_\mu$ be the unique box with $t(i,j)=h$. Define the \emph{content} of that label by
	\[ \operatorname{content}_t(h):=j-i. \]
\end{definition}

\begin{example}\label{ex:tableau-contents}
	In the following pair of diagrams, the left-hand diagram gives the labels of a standard tableau $t$, and the right-hand diagram gives the content of the label in each box:
	\[
		t=\ytableaushort{125,36,4} \qquad\longmapsto\qquad \begin{ytableau}
			0&1&2\\
			-1&0\\
			-2
		\end{ytableau}\;.
	\]
	Consequently,
	\[
		\begin{array}{c|rrrrrr}
			h                           & 1 & 2 & 3  & 4  & 5 & 6 \\\hline
			\operatorname{content}_t(h) & 0 & 1 & -1 & -2 & 2 & 0
		\end{array}
	\]
	For instance, the labels $4$, $5$, and $6$ occupy the boxes $(3,1)$, $(1,3)$, and $(2,2)$, respectively, so
	\[ \operatorname{content}_t(4)=1-3=-2,\qquad \operatorname{content}_t(5)=3-1=2,\qquad \operatorname{content}_t(6)=2-2=0. \]
\end{example}
Moving one box to the right increases its content by one; moving one box down decreases it by one. Boxes on a diagonal parallel to the main diagonal have the same content, and the main diagonal has content zero.

\begin{example}\label{ex:square-tableau-contents}
	For the square tableaux, the contents of the labels $2$ and $3$ are interchanged:
	\[
		\renewcommand{\arraystretch}{1.3} \begin{array}{c|cc}
		&\hspace{3pt}\raisebox{0pt}[\dimexpr\height+4pt\relax][\dimexpr\depth+4pt\relax]{$t_1=\ytableaushort{12,34}$}\hspace{3pt}&\hspace{3pt}\raisebox{0pt}[\dimexpr\height+4pt\relax][\dimexpr\depth+4pt\relax]{$t_2=\ytableaushort{13,24}$}\hspace{3pt}\\[7pt]\hline
		\text{content of }2&1&-1\\
		\text{content of }3&-1&1\\
		\text{content of }3\text{ minus content of }2&-2&2
		\end{array}
	\]
	Each column uses the tableau named in its heading. In particular, for the adjacent swap $(2\;3)t_2=t_1$, the content difference $\operatorname{content}_{t_2}(3)-\operatorname{content}_{t_2}(2)$ is $2$. This is the denominator appearing in Young's seminormal formula in the proof of \Cref{prop:seminormal}.
\end{example}

Fix $(a,k)\in\mathbb{Z}_{\ge1}^2$, put $d:=ak$, and define $\lambda:=(k,\ldots,k)\in\mathbb N^a$.
\begin{definition}\label{def:monic-seminormal-vectors}
	Choose a total order $\preceq$ on $\operatorname{SYT}(\lambda)$ such that
	\[ \forall\,(t,u)\in\operatorname{SYT}(\lambda)^2,\quad t\unrhd u\ \Longrightarrow\ u\preceq t. \]
	Define its \emph{strict part} $\prec$ by
	\[ \forall\,(u,t)\in\operatorname{SYT}(\lambda)^2,\quad u\prec t\ \Longleftrightarrow (u\preceq t\ \text{and}\ u\ne t). \]
	Thus tableaux with less dominant prefix shapes come first, and the sums below use the strict relation $\prec$. The least element for $\preceq$ is the column-reading tableau $\tau_0:=t_{\mathrm{col}}$ used in \eqref{eq:column-tableau-determinant}. Using the inner product induced from $(\mathbb{R}^a)^{\otimes d}$, define vectors in $L_{a,d}\otimes_\mathbb{Z}\mathbb{Q}$ recursively by
	\begin{equation}\label{eq:monic-gram-schmidt}
		\begin{gathered}
			s_{\tau_0}:=v_{\tau_0},\\
			\forall\,t\in\operatorname{SYT}(\lambda)\setminus\{\tau_0\},\quad
			s_t:=v_t-\sum_{\substack{u\in\operatorname{SYT}(\lambda)\\u\prec t}}
			\frac{\langle v_t,s_u\rangle}{\|s_u\|^2}\,s_u.
		\end{gathered}
	\end{equation}
	This is Gram--Schmidt orthogonalization: at each step, we subtract from $v_t$ its orthogonal projection onto the span of the preceding vectors. Since the standard Specht vectors form a basis, the resulting vectors are nonzero and mutually orthogonal. The recursion also gives
	\[ s_t-v_t\in\operatorname{span}_{\mathbb{Q}} \{v_u:u\in\operatorname{SYT}(\lambda),\ u\prec t\}. \]
	In particular, the coefficient of $v_t$ in the expansion of $s_t$ in the standard basis is one. We call $s_t$ the \emph{monic seminormal vector} indexed by $t$, where \emph{monic} refers to this coefficient. The resulting basis is Young's seminormal basis in the normalization of \cite[Section 2, the Gram--Schmidt formula preceding equation (1) and the following paragraph]{rh}; that paragraph also states its independence of the chosen total order $\preceq$ satisfying the displayed compatibility condition.
\end{definition}
\begin{definition}\label{def:seminormal-row-gaps}
	Let $\tau\in\operatorname{SYT}(\lambda)$. For $j\in\{0,\ldots,d\}$, define $\lambda^{(j)}:=\operatorname{sh}_j(\tau)$ and use its zero extension when taking row coordinates. For $j\in\{1,\ldots,d\}$, let $r_j$ be the row containing $j$ and define, using the preceding shape,
	\[ \forall\,i\in\{1,\ldots,r_j-1\},\quad \delta_{i,r_j}(j):= \lambda_i^{(j-1)}-\lambda_{r_j}^{(j-1)}+r_j-i. \]
	A legal addition implies $\delta_{i,r_j}(j)\ge2$.
\end{definition}

\begin{definition}\label{def:uniform-tableau-measure}
	Equip $\operatorname{SYT}(\lambda)$ with its power-set $\sigma$-algebra, and define the \emph{uniform probability measure} $\operatorname{\mathbb{P}}$ by
	\[ \forall\,\tau\in\operatorname{SYT}(\lambda),\quad \operatorname{\mathbb{P}}(\{\tau\}):=\frac1{\#(\operatorname{SYT}(\lambda))}. \]
	Write $\operatorname{\mathbb{E}}_{\operatorname{\mathbb{P}}}$ for expectation with respect to $\operatorname{\mathbb{P}}$.
\end{definition}

\begin{proposition}\label{prop:seminormal}
	For every $\tau\in\operatorname{SYT}(\lambda)$, let $r_j$ and $\delta_{i,r_j}(j)$ be the row indices and gaps of \Cref{def:seminormal-row-gaps}. The vector $s_\tau$ defined by \eqref{eq:monic-gram-schmidt} satisfies
	\begin{equation}\label{eq:pathnorm}
		\|s_\tau\|^2 =\prod_{j=1}^d\prod_{i=1}^{r_j-1} \frac{\delta_{i,r_j}(j)}{\delta_{i,r_j}(j)-1}.
	\end{equation}
	Moreover, the Gram matrix $G_{\mathrm{std}}$ of \Cref{def:standard-specht-gram} satisfies
	\begin{equation}\label{eq:grammean}
		\frac{\ln(\operatorname{det}({G_{\mathrm{std}}}))}{\#(\operatorname{SYT}(\lambda))} =\operatorname{\mathbb{E}}_{\operatorname{\mathbb{P}}}[\ln(\|s_\tau\|^2)].
	\end{equation}
\end{proposition}
\begin{proof}
	We compare the squared norm with the product on the right-hand side of \eqref{eq:pathnorm}. For each $\tau\in\operatorname{SYT}(\lambda)$, define
	\[ R(\tau):=\prod_{j=1}^d\prod_{i=1}^{r_j-1} \frac{\delta_{i,r_j}(j)}{\delta_{i,r_j}(j)-1}, \]
	where the row indices and gaps are computed from $\tau$. We take an empty product to be one. Every factor is positive, since $\delta_{i,r_j}(j)\ge2$. We will verify the equality $R(\tau)=\|s_\tau\|^2$ first at the column-reading tableau and then show that both sides change by the same ratio under every adjacent swap that preserves standardness.

	Let $\tau_0$ be the column-reading tableau of \Cref{def:monic-seminormal-vectors}; explicitly,
	\[ \forall\,(r,c)\in\{1,\ldots,a\}\times\{1,\ldots,k\},\quad \tau_0(r,c)=a(c-1)+r. \]
	Fix such a box $(r,c)$ and put $j:=a(c-1)+r$. Immediately before adding this box, the first $r-1$ rows contain $c$ boxes each, and the remaining rows contain $c-1$ boxes each. Thus the prefix shapes computed from $\tau_0$ satisfy
	\[
		\forall\,i\in\{1,\ldots,a\},\quad \lambda_i^{(j-1)}= \begin{cases}
			c,   & i\in\{1,\ldots,r-1\}, \\
			c-1, & i\in\{r,\ldots,a\}.
		\end{cases}
	\]
	Consequently,
	\[ \forall\,i\in\{1,\ldots,r-1\},\quad \delta_{i,r}(j)=c-(c-1)+r-i=r-i+1, \]
	and the factors contributed by this box have product
	\[ \prod_{i=1}^{r-1}\frac{\delta_{i,r}(j)}{\delta_{i,r}(j)-1} =\prod_{i=1}^{r-1}\frac{r-i+1}{r-i}=r. \]
	Multiplying over the $a$ boxes in each of the $k$ columns gives
	\[ R(\tau_0)=\prod_{c=1}^k\prod_{r=1}^a r=(a!)^k. \]
	On the other hand, \eqref{eq:monic-gram-schmidt} and \eqref{eq:column-tableau-determinant} give $s_{\tau_0}=v_{\tau_0}=\Delta_a^{\otimes k}$. The tensor $\Delta_a$ is a signed sum of $a!$ distinct orthonormal tensor words, so $\|\Delta_a\|^2=a!$. Hence
	\[ \|s_{\tau_0}\|^2=\|\Delta_a\|^{2k}=(a!)^k=R(\tau_0). \]

	We next compare two standard tableaux related by an adjacent swap. Let $j\in\{1,\ldots,d-1\}$ and suppose that $(\tau,\sigma)\in\operatorname{SYT}(\lambda)^2$ satisfy $\sigma=(j\;j+1)\tau$. The two labels cannot lie in the same row, since swapping them would violate its increasing order. After interchanging the names of the two tableaux if necessary, let $r\in\{2,\ldots,a\}$ be the row containing $j$ in $\tau$, and let $q\in\{1,\ldots,r-1\}$ be the row containing $j+1$. Put $\nu:=\operatorname{sh}_{j-1}(\tau) =\operatorname{sh}_{j-1}(\sigma)$. The boxes in question have coordinates $(r,\nu_r+1)$ and $(q,\nu_q+1)$, respectively. The prefix shapes of $\tau$ and $\sigma$ agree at every level other than $j$. At level $j$, one has
	\[
		\forall\,m\in\mathbb{Z}_{\ge1},\quad \sum_{i=1}^m\operatorname{sh}_j(\sigma)_i -\sum_{i=1}^m\operatorname{sh}_j(\tau)_i =\begin{cases}
			1, & m\in\{q,\ldots,r-1\}, \\
			0, & \text{otherwise}.
		\end{cases}
	\]
	Thus $\sigma\unrhd\tau$ and $\sigma\ne\tau$.

	Define $\Delta:=\nu_q-\nu_r+r-q$. Adding the box in row $r$ to $\nu$ preserves the partition condition, so $\nu_{r-1}\ge\nu_r+1$. Since $q<r$, one has $\nu_q\ge\nu_{r-1}$ and therefore $\Delta\ge r-q+1\ge2$. The contents in $\tau$ give
	\[ \operatorname{content}_\tau(j+1) -\operatorname{content}_\tau(j) =(\nu_q+1-q)-(\nu_r+1-r)=\Delta. \]
	The isometric identification in \Cref{rem:standard-specht-normalization} sends our vectors $v_t$ to the standard polytabloids used in \cite[Section 2, pp. 3--4]{rh}. It therefore also identifies the vectors obtained by the same monic Gram--Schmidt recursion. Young's seminormal formula (cf. \cite[Section 2, equation (1), p. 4]{rh}), in the case where the swapped tableau dominates the original one, gives
	\[ (j\;j+1)s_\tau=\rho s_\tau+s_\sigma, \qquad \rho:=\frac1\Delta. \]
	The transposition acts isometrically, and the two seminormal vectors are orthogonal. Taking squared norms therefore yields
	\[ \|s_\tau\|^2=\rho^2\|s_\tau\|^2+\|s_\sigma\|^2, \]
	and hence
	\begin{equation}\label{eq:normratio}
		\frac{\|s_\sigma\|^2}{\|s_\tau\|^2} =1-\rho^2=1-\frac1{\Delta^2}.
	\end{equation}

	We now compute the ratio $R(\sigma)/R(\tau)$. All additions before these two boxes are identical, and after both boxes have been added the shapes agree again. Their factors therefore cancel in the ratio. In $\tau$, the box in row $r$ is added first, whereas in $\sigma$ the box in row $q$ is added first. The factors for the addition in row $q$ depend only on rows with indices at most $q$; adding a box in row $r>q$ does not change any of them. For the addition in row $r$, the only changed upper row is $q$. Thus the two contributions from row $r$ are, respectively,
	\[
		\begin{aligned}
			\tau:\quad   &
			\frac{\Delta}{\Delta-1}
			\prod_{\substack{i\in\{1,\ldots,r-1\}\\i\ne q}}
			\frac{\nu_i-\nu_r+r-i}{\nu_i-\nu_r+r-i-1}, \\
			\sigma:\quad &
			\frac{\Delta+1}{\Delta}
			\prod_{\substack{i\in\{1,\ldots,r-1\}\\i\ne q}}
			\frac{\nu_i-\nu_r+r-i}{\nu_i-\nu_r+r-i-1}.
		\end{aligned}
	\]
	Cancellation gives
	\[ \frac{R(\sigma)}{R(\tau)} =\frac{(\Delta+1)/\Delta}{\Delta/(\Delta-1)} =\frac{\Delta^2-1}{\Delta^2} =\frac{\|s_\sigma\|^2}{\|s_\tau\|^2}. \]
	This equality also gives the corresponding ratios for the reverse swap by taking reciprocals.

	It remains to explain why these comparisons reach every standard tableau. On the set $D_\lambda$ of boxes, define the partial order
	\[ \forall\,((i,c),(r,b))\in D_\lambda^2,\quad (i,c)\le_D(r,b)\ \Longleftrightarrow\ i\le r\text{ and }c\le b. \]
	A standard tableau lists the boxes in an order extending $\le_D$. Indeed, distinct boxes comparable under $\le_D$ can be joined, in the direction of that order, by steps to the right and down inside $D_\lambda$. The row and column inequalities therefore put the smaller box first. Conversely, labelling the boxes in any order extending $\le_D$ gives a standard tableau. We show directly that any two such orders can be joined by swapping consecutive incomparable boxes. Let $(x_1,\ldots,x_d)$ be the current order and $(y_1,\ldots,y_d)$ the desired order. Suppose that their first $m-1$ terms already agree, where $m\in\{1,\ldots,d\}$, and let $h\in\{m,\ldots,d\}$ satisfy $x_h=y_m$. For every $p\in\{m,\ldots,h-1\}$, the box $x_p$ is incomparable with $y_m$. Indeed, $y_m\le_D x_p$ would contradict the current order, while $x_p\le_D y_m$ would contradict the desired order, since $x_p$ is not among its first $m-1$ terms. Moving $y_m$ successively to position $m$ therefore uses only swaps of consecutive incomparable boxes. Each such swap preserves the property of extending $\le_D$. Repeating for $m\in\{1,\ldots,d\}$ gives the desired order. In tableau notation these are exactly adjacent swaps of labels that preserve standardness. Along every such swap, $\|s_\tau\|^2/R(\tau)$ is unchanged. The initial calculation at $\tau_0$ now gives
	\[ \forall\,\tau\in\operatorname{SYT}(\lambda),\quad \frac{\|s_\tau\|^2}{R(\tau)} =\frac{\|s_{\tau_0}\|^2}{R(\tau_0)}=1, \]
	which proves \eqref{eq:pathnorm}.

	Finally, permuting a basis does not change its Gram determinant, so we may order the standard Specht basis by the total order $\preceq$ of \Cref{def:monic-seminormal-vectors}. Let $C$ be the matrix whose column indexed by $t$ consists of the coefficients of $s_t$ in this ordered basis. The recursion \eqref{eq:monic-gram-schmidt} makes $C$ triangular with every diagonal entry equal to one, so $\operatorname{det}(C)=1$. The change of basis rule for Gram matrices and the orthogonality of the seminormal basis give
	\[ C^{\mathsf T}G_{\mathrm{std}}C =\operatorname{diag}\bigl((\|s_t\|^2)_{t\in\operatorname{SYT}(\lambda)}\bigr). \]
	Taking determinants yields
	\[ \operatorname{det}({G_{\mathrm{std}}}) =\prod_{\tau\in\operatorname{SYT}(\lambda)}\|s_\tau\|^2. \]
	Taking logarithms and dividing by the number of standard tableaux gives
	\[ \frac{\ln(\operatorname{det}(G_{\mathrm{std}}))}{\#(\operatorname{SYT}(\lambda))} =\frac1{\#(\operatorname{SYT}(\lambda))} \sum_{\tau\in\operatorname{SYT}(\lambda)}\ln(\|s_\tau\|^2). \]
	By \Cref{def:uniform-tableau-measure}, the right-hand side is precisely $\operatorname{\mathbb{E}}_{\operatorname{\mathbb{P}}}[\ln(\|s_\tau\|^2)]$. This proves \eqref{eq:grammean}.
\end{proof}

\begin{example}\label{ex:square-seminormal}
	Take $(a,k)=(2,2)$ and the tableaux in \eqref{eq:tableaux-square}. Then $\tau_0=t_2\prec t_1$ and $\sigma=t_1=(2\;3)t_2$. The inner products in \Cref{ex:square-specht-lattice} give the explicit Gram--Schmidt step
	\[ s_{t_2}=v_{t_2},\qquad s_{t_1}=v_{t_1} -\frac{\langle v_{t_1},v_{t_2}\rangle}{\|v_{t_2}\|^2}v_{t_2} =v_{t_1}-\frac12v_{t_2}. \]
	Their content difference is $2$, so \eqref{eq:normratio} gives $\|s_{t_1}\|^2=(1-1/4)\|s_{t_2}\|^2=3$, since $\|s_{t_2}\|^2=\|v_{t_2}\|^2=4$. Thus the product of the two squared seminormal norms is $12$, in agreement with the Gram matrix computed in \Cref{subsec:specht-tensors}.
\end{example}

\section{Sublinear covolume from exact prefix probabilities}\label{sec:asymptotic}

The main estimate of this section is the following.

\begin{proposition}\label{prop:sublinear}
	For every fixed $a\in\mathbb{Z}_{\ge1}$,
	\begin{equation}\label{eq:sublinear}
		\lim\limits_{\substack{d\to\infty\\d\in a\mathbb{Z}_{\ge1}}} \frac{-\operatorname{\widehat{\mu}}(\overline L_{a,d}^{\,S})}{d}=0.
	\end{equation}
\end{proposition}

Fix $(a,k)\in\mathbb{Z}_{\ge2}\times\mathbb{Z}_{\ge1}$, put $d:=ak$, and define $\lambda:=(k,\ldots,k)\in\mathbb N^a$. We keep $a$ fixed when taking limits as $k$ tends to infinity. Let $\tau$ be a random element of $\operatorname{SYT}(\lambda)$ with the uniform probability law $\operatorname{\mathbb{P}}$ of \Cref{def:uniform-tableau-measure}, and use its expectation $\operatorname{\mathbb{E}}_{\operatorname{\mathbb{P}}}$.

\subsection{From covolume to insertion costs} \label{subsec:insertion-costs}

To prove \Cref{prop:sublinear}, we first reduce the slope estimate to an average over standard tableaux. By \eqref{eq:covol-gram-bound} and \eqref{eq:grammean},
\begin{equation}\label{eq:average-covolume-reduction}
	0\le\frac{2\ln(\operatorname{covol}(L_{a,d}))}{d\#(\operatorname{SYT}(\lambda))} \le\frac{\operatorname{\mathbb{E}}_{\operatorname{\mathbb{P}}}[\ln(\|s_\tau\|^2)]}{d}.
\end{equation}
By \eqref{eq:lattice-extension-slope}, the negative slope $-\operatorname{\widehat{\mu}}(\overline L_{a,d}^{\,S})$ is $\nu(\Omega_\infty)$ times the logarithmic covolume divided by $\#(\operatorname{SYT}(\lambda))$. It therefore suffices to prove
\begin{equation}\label{eq:average-cost-target}
	\lim\limits_{k\to\infty} \frac{\operatorname{\mathbb{E}}_{\operatorname{\mathbb{P}}}[\ln(\|s_\tau\|^2)]}{ak}=0.
\end{equation}
The averaging in this statement is essential. For the column-reading tableau $\tau_0$, \eqref{eq:column-tableau-norm} and $s_{\tau_0}=v_{\tau_0}$ give
\[ \|s_{\tau_0}\|^2=(a!)^k,\qquad \frac{\ln(\|s_{\tau_0}\|^2)}{ak} =\frac{\ln(a!)}a>0. \]
Thus a bound obtained by taking the largest seminormal norm cannot prove \eqref{eq:average-cost-target}.

A standard tableau $\tau$ determines a growth path
\[ \varnothing=\lambda^{(0)}\subsetneq\lambda^{(1)} \subsetneq\cdots\subsetneq\lambda^{(d)}=\lambda, \qquad \lambda^{(t)}:=\operatorname{sh}_t(\tau), \]
where the $j$-th step adds the box labelled $j$. For $(\tau,j)\in\operatorname{SYT}(\lambda)\times\{1,\ldots,d\}$, let $r_j$ be the row of this box. Using the quantities $\delta_{i,r_j}(j)$ from \Cref{sec:seminormal}, define its \emph{logarithmic insertion cost} by
\begin{equation}\label{eq:insertion-cost}
	H_j(\tau):=\sum_{i=1}^{r_j-1} \ln\left(\frac{\delta_{i,r_j}(j)} {\delta_{i,r_j}(j)-1}\right).
\end{equation}
Equation \eqref{eq:pathnorm} then gives
\begin{equation}\label{eq:total-insertion-cost}
	\ln(\|s_\tau\|^2)=\sum_{j=1}^d H_j(\tau),\qquad \frac{\operatorname{\mathbb{E}}_{\operatorname{\mathbb{P}}}[\ln(\|s_\tau\|^2)]}{d} =\frac1d\sum_{j=1}^d\operatorname{\mathbb{E}}_{\operatorname{\mathbb{P}}}[H_j(\tau)].
\end{equation}
The problem is therefore to show that the average cost per insertion tends to zero.

For example, suppose that $a=2$ and the current prefix has shape $(u,v)\in\mathbb{Z}_{\ge0}^2$, with $k\ge u\ge v$. Put $h:=u-v$. An insertion in the first row has cost zero. An insertion in the second row is legal only when $h\ge1$, and its cost is
\[ \ln\left(\frac{h+1}{h}\right). \]
Thus this cost is small when the two row lengths are far apart. For $k=2$, the tableau $t_1=\ytableaushort{12,34}$ gives the path
\[ \varnothing\xrightarrow{\ 0\ }(1) \xrightarrow{\ 0\ }(2) \xrightarrow{\ \ln(3/2)\ }(2,1) \xrightarrow{\ \ln(2)\ }(2,2), \]
where each arrow is labelled by its insertion cost. Their sum is $\ln(3)$, in agreement with $\|s_{t_1}\|^2=3$ computed in \Cref{prop:seminormal}.

For general $a$, define the shifted prefix coordinates by
\[ \forall\,t\in\{0,\ldots,d\},\quad b^{(t)}:=(\lambda_i^{(t)}+a-i)_{i=1}^a. \]
They satisfy $b_1^{(t)}>\cdots>b_a^{(t)}$. At time $t=j-1$,
\begin{equation}\label{eq:cost-prefix-gaps}
	\forall\,i\in\{1,\ldots,r_j-1\},\quad \delta_{i,r_j}(j)=b_i^{(t)}-b_{r_j}^{(t)}.
\end{equation}
Consequently, for every $L\in[2,\infty)$,
\begin{equation}\label{eq:insertion-cost-bounds}
	\begin{aligned}
		0\le H_j(\tau) & \le(a-1)\ln(2),                                 \\
		\min_{\substack{(i,r)\in\{1,\ldots,a\}^2\\i<r}}|b_i^{(j-1)}-b_r^{(j-1)}|>L
		               & \ \Longrightarrow\ H_j(\tau)\le\frac{a-1}{L-1}.
	\end{aligned}
\end{equation}
The first bound follows from $\delta_{i,r_j}(j)\ge2$; the second follows from $\ln(1+x)\le x$ for $x\in\mathbb{R}_{\ge0}$. It holds for every legal choice of the next row. Hence the prefix distribution alone will suffice to control the cost.

We shall prove that, for each fixed $(\varepsilon,L)\in(0,1/2)\times[2,\infty)$, the probability of a pairwise gap at most $L$ tends to zero uniformly over $t\in\mathbb{Z}\cap[\varepsilon d,(1-\varepsilon)d]$. This restriction to intermediate times is necessary: at both $t=0$ and $t=d$, consecutive shifted coordinates differ by one. The initial and final time intervals will instead be controlled by their small proportion of the total number of steps.

\subsection{The exact distribution of a prefix} \label{subsec:prefix-distribution}

To estimate the gaps, we first count tableaux with a prescribed prefix shape. For a shape $\nu\subseteq\lambda$, the prefix and its completion can be filled independently. The rectangular outer shape allows us to turn the completion into an ordinary Young diagram by rotation and reversal of labels. This gives the product of tableau counts in \eqref{eq:prefix} below. Applying the shifted hook-length formula to both counts will produce one Vandermonde factor from each, and hence a \emph{squared Vandermonde weight}.

For this calculation, define
\[ c:=\frac{a(a-1)}2,\qquad M:=k+a-1,\qquad V_0:=\prod_{j=1}^{a-1}j!. \]
For $t\in\{0,\ldots,d\}$, put $A:=t+c$ and define
\[ \mathcal X_{k,t}:= \left\{b\in\{0,\ldots,M\}^a:\sum_{i=1}^a b_i=A\right\}. \]
For $b\in\mathcal X_{k,t}$, define its Vandermonde product by
\[ V(b):=\prod_{\substack{(i,j)\in\{1,\ldots,a\}^2\\i<j}}(b_i-b_j). \]
Equip $\mathcal X_{k,t}$ with its power-set $\sigma$-algebra. Define a probability mass function on it by
\[ B_{k,t}(b):=\frac{\prod_{i=1}^a\binom M{b_i}}{\binom{aM}A}. \]
It is the \emph{multivariate hypergeometric distribution} obtained by choosing $A$ objects uniformly from $a$ groups of $M$ objects and recording the number chosen from each group. In particular, its masses sum to one. Write $\operatorname{\mathbb{E}}_{B_{k,t}}$ and $\operatorname{\mathbb{P}}_{B_{k,t}}$ for expectation and probability under this distribution. Define
\[ Z_{k,t}:=\operatorname{\mathbb{E}}_{B_{k,t}}[V(b)^2]. \]

\begin{samepage}
	\begin{proposition}\label{prop:prefix-law}
		Let $t\in\{0,\ldots,d\}$, let $\nu\subseteq\lambda$ be a partition of $t$, and use its zero extension to write $(\nu_1,\ldots,\nu_a)$. Since $\lambda=(k,\ldots,k)$ has $a$ parts, this containment means that $\nu$ has at most $a$ positive parts and $\nu_i\in\{0,\ldots,k\}$ for every $i\in\{1,\ldots,a\}$. Define
		\[ \nu^c:=(k-\nu_a,\ldots,k-\nu_1),\qquad \forall\,i\in\{1,\ldots,a\},\quad b_i:=\nu_i+a-i. \]
		The prefix shape $\lambda^{(t)}:=\operatorname{sh}_t(\tau)$ satisfies
		\begin{equation}\label{eq:prefix}
			\operatorname{\mathbb{P}}(\lambda^{(t)}=\nu)=\frac{\#(\operatorname{SYT}(\nu))\cdot \#(\operatorname{SYT}(\nu^c))}{\#(\operatorname{SYT}(\lambda))}.
		\end{equation}
		Moreover,
		\begin{align}
			Z_{k,t}
			 & =a!V_0\frac{(M!)^a}{\prod_{j=0}^{a-1}(k+j)!}
			\cdot\frac{\binom dt}{\binom{aM}{t+c}}\notag    \\
			 & =a!V_0\prod_{j=1}^{a-1}(k+j)^j
			\frac{\binom dt}{\binom{d+2c}{t+c}}.
			\label{eq:normalizer}
		\end{align}
		If a tuple $b\in\mathcal X_{k,t}$ is sampled with masses $V(b)^2B_{k,t}(b)/Z_{k,t}$ and then sorted in decreasing order, its distribution is that of $(\lambda_i^{(t)}+a-i)_{i=1}^a$.
	\end{proposition}
\end{samepage}
\begin{proof}
	Fix a partition $\nu\subseteq\lambda$ of $t$. There are $\#(\operatorname{SYT}(\nu))$ standard fillings of its boxes by the labels $\{1,\ldots,t\}$. The remaining boxes form the skew diagram $D_\lambda\setminus D_\nu$. Rotation by $180$ degrees sends a box $(i,j)$ to $(a+1-i,k+1-j)$ and sends this skew diagram to $D_{\nu^c}$. Replace a label $h\in\{t+1,\ldots,d\}$ by $d+1-h$. The rotation reverses row and column order, and the replacement reverses label order, so the resulting filling is a standard tableau of shape $\nu^c$. The operation is invertible. There are therefore $\#(\operatorname{SYT}(\nu^c))$ choices for the remaining labels. All labels in $D_\nu$ are smaller than those in its complement; since $D_\nu$ is a Young diagram, there are no further boundary conditions. Dividing the resulting count by $\#(\operatorname{SYT}(\lambda))$ proves \eqref{eq:prefix}, including $t=0$ and $t=d$ under the empty-partition convention of \Cref{sec:seminormal}.

	The shifted coordinates satisfy
	\[ b\in\mathcal X_{k,t},\qquad b_1>\cdots>b_a. \]
	Conversely, a strictly decreasing tuple in $\mathcal X_{k,t}$ gives $\nu_i:=b_i-a+i$. Its coordinates are weakly decreasing. Since the tuple has $a$ distinct entries in $\{0,\ldots,k+a-1\}$, one has $a-i\le b_i\le k+a-i$, so $\nu_i\in\{0,\ldots,k\}$ and $\nu\subseteq\lambda$. This establishes the claimed correspondence of supports.

	Applying \eqref{eq:shifted-count} to the prefix partition gives
	\[ \#(\operatorname{SYT}(\nu))=\frac{t!V(b)}{\prod_{i=1}^a b_i!}. \]

	The shifted coordinates of $\nu^c$ are $(M-b_{a+1-i})_{i=1}^a$. Reversing the indices in the product of positive coordinate differences preserves the Vandermonde product, and hence
	\[ \#(\operatorname{SYT}(\nu^c))=\frac{(d-t)!V(b)}{\prod_{i=1}^a(M-b_i)!},\qquad \#(\operatorname{SYT}(\lambda))=\frac{d!V_0}{\prod_{j=0}^{a-1}(k+j)!}. \]
	Substitution in \eqref{eq:prefix} yields, for a strictly decreasing $b\in\mathcal X_{k,t}$,
	\[ \operatorname{\mathbb{P}}\bigl((\lambda_i^{(t)}+a-i)_{i=1}^a=b\bigr) =C_{k,t}V(b)^2\prod_{i=1}^a\binom M{b_i},\qquad C_{k,t}:=\frac{\prod_{j=0}^{a-1}(k+j)!} {\binom dt\,V_0(M!)^a}. \]
	This is the exact relation between tableau counting and sampling without replacement.

	Both $B_{k,t}(b)$ and $V(b)^2$ are symmetric in the coordinates. A tuple with a repeated coordinate has zero Vandermonde weight; a tuple with distinct coordinates has exactly $a!$ different permutations. Summing the last probability identity over the strictly decreasing tuples therefore gives
	\[ 1=C_{k,t}\binom{aM}A\,\frac{Z_{k,t}}{a!}. \]
	This proves the first equality in \eqref{eq:normalizer} and the sorting assertion. The second equality follows from $aM=d+2c$ and cancellation of factorials:
	\[ \frac{(M!)^a}{\prod_{j=0}^{a-1}(k+j)!} =\prod_{r=0}^{a-1}\prod_{j=r+1}^{a-1}(k+j) =\prod_{j=1}^{a-1}(k+j)^j. \]
\end{proof}

For any event $\mathcal E\subseteq\mathcal X_{k,t}$ invariant under coordinate permutations, \Cref{prop:prefix-law} gives the useful identity
\begin{equation}\label{eq:weighted-event}
	\operatorname{\mathbb{P}}\bigl((\lambda_i^{(t)}+a-i)_{i=1}^a\in\mathcal E\bigr) =\frac{\operatorname{\mathbb{E}}_{B_{k,t}}[\mathbf1_{\mathcal E}(b)V(b)^2]}{Z_{k,t}}.
\end{equation}
Here $\mathbf1_{\mathcal E}$ denotes the indicator function. The formula separates two parts of the prefix law: the hypergeometric distribution controls the tails, and the weight $V(b)^2$ suppresses small coordinate gaps. The normalization $Z_{k,t}$ is essential for turning these observations into a probability bound. In particular, a bound for the unweighted hypergeometric law must be transferred through \eqref{eq:weighted-event} before it can be used for tableaux.

\subsection{Small gaps at intermediate times} \label{subsec:prefix-small-gaps}

The next goal is a uniform estimate over the intermediate prefix times. Uniformity matters because the number of terms in \eqref{eq:total-insertion-cost} grows with $k$. For $(\varepsilon,L)\in(0,1/2)\times\mathbb{R}_{>0}$, define
\[ \mathcal T_{k,\varepsilon}:= \mathbb{Z}\cap[\varepsilon ak,(1-\varepsilon)ak]. \]
Define
\[ p_{k,\varepsilon,L}:= \max\left(\{0\}\cup \left\{\operatorname{\mathbb{P}}\left( \min_{\substack{(i,j)\in\{1,\ldots,a\}^2\\i<j}} |b_i^{(t)}-b_j^{(t)}|\le L\right): t\in\mathcal T_{k,\varepsilon}\right\}\right). \]
Including zero specifies the value when the set of intermediate integer times is empty.

The proof of \Cref{lem:smallgap} will split the weighted probability into a central region and its complement. In the central region, all coordinate differences have a controlled size; requiring one of them to be at most $L$ makes one factor of $V(b)^2$ small. On the complement, the polynomial bound on $V(b)^2$ is dominated by an exponential tail estimate for the hypergeometric distribution. We first establish that estimate.

\begin{lemma}\label{lem:sampling-tail}
	Let $(N,M,A)\in\mathbb{Z}_{\ge1}\times\mathbb{Z}_{\ge0}\times\mathbb{Z}_{\ge1}$. Assume that $M\le N$ and $A\le N$. Choose an $A$-element subset uniformly from a population of size $N$, of which $M$ elements are marked. Let $X$ be the number of marked elements chosen, let $B$ denote its probability law on $\{0,\ldots,A\}$ with the power-set $\sigma$-algebra, and put $p:=M/N$. For every $u\in\mathbb{R}_{>0}$,
	\[ \operatorname{\mathbb{P}}_B(|X-Ap|\ge u)\le 2\operatorname{exp}\left(-\frac{u^2}{4A}\right). \]
\end{lemma}
\begin{proof}
	For $(x,r)\in\mathbb{Z}_{\ge0}^2$, define $(x)_r:=x(x-1)\cdots(x-r+1)$, with $(x)_0:=1$. Counting pairs consisting of a sample and an $r$-element marked subset of that sample gives, for $r\in\{0,\ldots,A\}$,
	\[ \operatorname{\mathbb{E}}_B\left[\binom Xr\right] =\binom Ar\frac{(M)_r}{(N)_r}\le\binom Ar p^r. \]
	For $r\le M$, the inequality follows from $(M-j)/(N-j)\le M/N$ for $j\in\{0,\ldots,r-1\}$; for $r>M$ the numerator vanishes. For $\theta\in\mathbb{R}_{\ge0}$, put $z:=e^\theta-1$. Expanding $(1+z)^X$ therefore gives
	\[ \forall\,\theta\in\mathbb{R}_{\ge0},\quad \operatorname{\mathbb{E}}_B[\operatorname{exp}(\theta X)]\le(1-p+pe^\theta)^A. \]
	Applying this inequality to the unmarked count $A-X$ gives the same bound for $\theta<0$. For $|\theta|\le1$, the elementary exponential remainder bound $e^\theta-1-\theta\le\theta^2$ implies
	\[ 1+p(e^\theta-1)\le1+p\theta+\theta^2 \le e^{p\theta+\theta^2}, \qquad \operatorname{\mathbb{E}}_B[\operatorname{exp}(\theta(X-Ap))]\le e^{A\theta^2}. \]
	If $0<u\le A$, Markov's inequality with $\theta=u/(2A)$ bounds the upper tail by $e^{-u^2/(4A)}$; the parameter $-\theta$ gives the same lower-tail bound. If $u>A$, the event is empty because $0\le X,Ap\le A$. This proves the assertion.
\end{proof}

\begin{lemma}\label{lem:smallgap}
	For each fixed $(\varepsilon,L)\in(0,1/2)\times\mathbb{R}_{>0}$,
	\[ \lim\limits_{k\to\infty}p_{k,\varepsilon,L}=0. \]
\end{lemma}
\begin{proof}
	Fix $k\in\mathbb{Z}_{\ge1}$ and $t\in\mathbb{Z}\cap[\varepsilon ak,(1-\varepsilon)ak]$. We first bound the denominator of \eqref{eq:weighted-event} from below. The exact normalization in \Cref{prop:prefix-law} is useful here because
	\[ \frac{\binom dt}{\binom{d+2c}{t+c}} =\frac{\prod_{j=1}^c(t+j)(d-t+j)} {\prod_{j=1}^{2c}(d+j)}. \]
	Since $d+2c=a(k+a-1)\le ad$ for $k\in\mathbb{Z}_{\ge1}$, \eqref{eq:normalizer} implies the uniform lower bound
	\begin{equation}\label{eq:Zlower}
		Z_{k,t}\ge a!V_0(\varepsilon/a)^{2c}k^c.
	\end{equation}

	Under $B_{k,t}$, each $b_i$ counts marked elements in a sample of size $A=t+c$ from $N=aM$ objects, with $M$ marked. Thus \Cref{lem:sampling-tail} gives
	\[ \forall\,(i,u)\in\{1,\ldots,a\}\times\mathbb{R}_{>0},\quad \operatorname{\mathbb{P}}_{B_{k,t}}(|b_i-A/a|\ge u)\le2e^{-u^2/(4A)}. \]
	We use a central region slightly wider than the square-root scale. Set $\eta:=1/(4c)$ and $u:=k^{1/2+\eta}$. This choice satisfies $-1+2\eta(c-1)<0$, which will make the central small-gap contribution tend to zero after division by $Z_{k,t}$. Since $A\le aM\le a^2k$, a union bound yields
	\begin{equation}\label{eq:tail}
		\operatorname{\mathbb{P}}_{B_{k,t}}\left(\max_{i\in\{1,\ldots,a\}}|b_i-A/a|\ge u\right) \le2a\operatorname{exp}\left(-\frac{k^{2\eta}}{4a^2}\right).
	\end{equation}

	We now transfer these estimates to the tableau law using \eqref{eq:weighted-event}. There are $c$ squared factors in $V(b)^2$. On the central event $\max_{i\in\{1,\ldots,a\}}|b_i-A/a|<u$, all pair differences are less than $2u$. If one difference is at most $L$, then
	\[ V(b)^2\le L^2(2u)^{2(c-1)}. \]
	The contribution of this event to the numerator in \eqref{eq:weighted-event} is at most this deterministic bound, since its $B_{k,t}$-probability is at most one. By \eqref{eq:Zlower}, its tableau probability is at most
	\[ \frac{2^{2c-2}a^{2c}}{a!V_0}\, \varepsilon^{-2c}L^2 k^{-(c+1)/(2c)}. \]
	Here we used $-1+2\eta(c-1)=-(c+1)/(2c)<0$. It remains to bound the complement of the central event. Here it suffices to use the global polynomial bound $V(b)^2\le M^{2c}\le(ak)^{2c}$ in the numerator of \eqref{eq:weighted-event}, and then apply \eqref{eq:Zlower} and \eqref{eq:tail}. Its tableau probability is at most
	\[ \frac{2a^{4c+1}}{a!V_0}\, \varepsilon^{-2c}k^c{\operatorname{exp}(-k^{1/(2c)}/(4a^2))}. \]
	Both upper bounds tend to zero for fixed $(\varepsilon,L)$ and do not depend on the chosen intermediate time $t$. Both events are invariant under sorting the coordinates. These bounds prove the asserted uniform limit for the prefix distribution.
\end{proof}

\subsection{Averaging the insertion costs} \label{subsec:average-insertion-costs}

We now apply \Cref{lem:smallgap} to \eqref{eq:total-insertion-cost}. The two endpoint time intervals contain few steps. Among the intermediate steps, a small gap has small probability, while the cost on its complement is small by \eqref{eq:insertion-cost-bounds}. These three contributions give the desired average estimate.

\begin{proof}[Proof of \Cref{prop:sublinear}]
	If $a=1$, the lattice $L_{1,d}=\mathbb{Z} e_1^{\otimes d}$ has rank one and its displayed basis vector has norm one. Hence $\operatorname{\widehat{\mu}}(\overline L_{1,d}^{\,S})=0$ for every $d\in\mathbb{Z}_{\ge1}$, proving the assertion in this case. Assume that $a\in\mathbb{Z}_{\ge2}$. Fix $(\varepsilon,L)\in(0,1/2)\times[2,\infty)$. Use the costs $H_j$ defined in \eqref{eq:insertion-cost} and split their expectations according to the prefix time $t=j-1$ and the gaps of $b^{(t)}$.
	\begin{enumerate}[label=\rm(\arabic*)]
		\item If $t\notin\mathcal T_{k,\varepsilon}$, use $H_j\le(a-1)\ln(2)$. There are at most $2\varepsilon d+2$ such times in $\{0,\ldots,d-1\}$.
		\item If $t\in\mathcal T_{k,\varepsilon}$ and some pairwise gap is at most $L$, use the same bound on $H_j$. By definition of $p_{k,\varepsilon,L}$, the contribution of this event to $\operatorname{\mathbb{E}}_{\operatorname{\mathbb{P}}}[H_j(\tau)]$ is at most $(a-1)\ln(2)\,p_{k,\varepsilon,L}$.
		\item If $t\in\mathcal T_{k,\varepsilon}$ and all pairwise gaps exceed $L$, then \eqref{eq:insertion-cost-bounds} gives $H_j\le(a-1)/(L-1)$. This holds for every legal next insertion, so the probability of its row need not be computed.
	\end{enumerate}
	Summing these contributions and using \eqref{eq:total-insertion-cost} gives
	\begin{equation}\label{eq:average-cost-bound}
		\begin{aligned}
			\frac1d\operatorname{\mathbb{E}}_{\operatorname{\mathbb{P}}}[\ln(\|s_\tau\|^2)]
			\le{} & (2\varepsilon+2/d)(a-1)\ln(2) \\
			      & +\frac{a-1}{L-1}
			+(a-1)\ln(2)\,p_{k,\varepsilon,L}.
		\end{aligned}
	\end{equation}
	The parameters $\varepsilon$ and $L$ are fixed while $k$ tends to infinity, so \Cref{lem:smallgap} applies to the last term. Taking the superior limit as $k$ tends to infinity yields
	\[ 0\le\limsup\limits_{k\to\infty} \frac{\operatorname{\mathbb{E}}_{\operatorname{\mathbb{P}}}[\ln(\|s_\tau\|^2)]}{ak} \le2\varepsilon(a-1)\ln(2)+\frac{a-1}{L-1}. \]
	Letting $L$ tend to infinity and then letting $\varepsilon$ tend to zero proves \eqref{eq:average-cost-target}, since the averaged costs are nonnegative. Finally, \eqref{eq:average-covolume-reduction} gives
	\[ \lim\limits_{k\to\infty} \frac{\ln(\operatorname{covol}(L_{a,ak}))}{ak\#(\operatorname{SYT}(\lambda))}=0. \]
	Multiplication by the fixed finite number $\nu(\Omega_\infty)$ yields \eqref{eq:sublinear}.
\end{proof}

\section{Adelic contractions and the geometric slope inequality} \label{sec:numberfields}
In this section and \Cref{sec:allrank}, assume $\operatorname{char}(K)=0$. The integral Specht lattices determine constant adelic bundles; their non-Archimedean norms are specified independently of their embeddings into tensor powers.

The distinction from the single determinant-tensor estimate in \cite[Proposition 5.3.1 and Remark 5.3.2]{cm} is the choice of auxiliary source. We retain the entire Specht representation, with its integral basis at non-Archimedean places and its induced tensor metric at Archimedean places. This representation is semistable under the action of its finite permutation group. The slope comparison therefore uses its average logarithmic covolume, rather than the norm of one selected determinant tensor. The following construction verifies that this choice is compatible with every local contraction used in the proof.

\begin{lemma}\label{lem:arithmetic-specht}
	Let $\overline A$ be a positive-rank Hermitian adelic bundle on $S$, put $a:=\operatorname{dim}_K(A)$, and fix $k\in\mathbb{Z}_{\ge1}$. Put $d:=ak$ and $\lambda:=(k,\ldots,k)\in\mathbb N^a$. Let $\overline L_{a,d}^{\,S}$ be the Hermitian adelic vector bundle obtained from the Euclidean lattice $L_{a,d}$ by extension of scalars, as in \Cref{subsec:finite-groups}. There is a canonical $K$-linear injection
	\[ {\iota_A:(\operatorname{det}(A))^{\otimes k} \otimes_K L_{a,d}^{\,S}\longrightarrow A^{\otimes d}} \]
	which is isometric at Archimedean places and contracting at all other places. Moreover,
	\begin{equation}\label{eq:arithmetic-specht-slope}
		{\operatorname{\widehat{\mu}}\left((\operatorname{det}(\overline A))^{\otimes k} \otimes\overline L_{a,d}^{\,S}\right)} =d\operatorname{\widehat{\mu}}(\overline A) +{\operatorname{\widehat{\mu}}(\overline L_{a,d}^{\,S})}.
	\end{equation}
	{The image of $\iota_A$ satisfies}
	\[ {\operatorname{Im}(\iota_A)=(A^{\otimes d})^{\operatorname{SL}(A)}}. \]
	{The permutation action of $\mathfrak S_d$ on $(\operatorname{det}(\overline A))^{\otimes k}\otimes\overline L_{a,d}^{\,S}$ is absolutely irreducible and isometric.}
\end{lemma}
\begin{proof}
	The construction in \Cref{subsec:finite-groups} shows that $\overline L_{a,d}^{\,S}$ is a Hermitian adelic vector bundle. {Its non-Archimedean norms are determined by the $\mathbb{Z}$-module structure of $L_{a,d}$ through \eqref{eq:scalar-extension-lattice-norm}.}

	For a basis $b=(b_1,\ldots,b_a)$ of $A$, define $\eta_b:=b_1\wedge\cdots\wedge b_a$. For $x\in{L_{a,d}^{\,S}}$, the rule
	\[ \eta_b^{\otimes k}\otimes x \longmapsto b^{\otimes d}x \]
	defines $\iota_A$. Here $b^{\otimes d}x$ is obtained by substituting $b_i$ for the $i$th tensor letter in $x$. By \eqref{eq:specht-determinant-orbit}, {the lattice is generated by permutations of products of $k$ determinant tensors, so} after scalar extension,
	\[ \forall\,g\in\operatorname{GL}_a(K),\quad g^{\otimes d}x=(\operatorname{det}(g))^k x. \]
	Consequently the rule is independent of $b$. The vectors $(v_t)_{t\in\operatorname{SYT}(\lambda)}$ are linearly independent over $\mathbb{Q}$ and remain so after extension of scalars to $K$. Thus $\iota_A$ is injective. At an Archimedean place an orthonormal basis makes the displayed rule isometric for the determinant metric and the canonical Specht metric. In particular, the norm $\sqrt{a!}$ of an unnormalized determinant tensor is retained in the Specht metric.

	Fix $\omega\in\Omega\setminus\Omega_\infty$ and $\alpha\in(0,1)$. Choose an $\alpha$-orthogonal basis $b=(b_1,\ldots,b_a)$ of $A_\omega$, so that
	\[
		\begin{gathered}
			\forall\,(z_1,\ldots,z_a)\in K_\omega^a,\\
			\alpha\max_{i\in\{1,\ldots,a\}}|z_i|_\omega\|b_i\|_\omega
			\le\left\|\sum_{i=1}^a z_i b_i\right\|_\omega
			\le\max_{i\in\{1,\ldots,a\}}|z_i|_\omega\|b_i\|_\omega.
		\end{gathered}
	\]
	Such bases exist for every finite-dimensional ultrametric normed space (cf. \cite[Section 1.2.3]{cm}). The determinant norm therefore satisfies
	\[ \|\eta_b\|_{\omega,\operatorname{det}} \ge\alpha^a\prod_{i=1}^a\|b_i\|_\omega. \]
	Every element of $L_{a,d}$ has integer coordinates in the basis of tensor words, supported on words using each letter exactly $k$ times. This follows from the same property of all the orbit generators. Writing $\|\mathord{\cdot}\|_{L,\omega}$ for the non-Archimedean Specht norm, {expansion in any integral basis of $L_{a,d}$ and \eqref{eq:scalar-extension-lattice-norm} give}, for every $x\in K_\omega\otimes_\mathbb{Z} L_{a,d}$,
	\[ \|b^{\otimes d}x\|_\omega \le\left(\prod_{i=1}^a\|b_i\|_\omega\right)^k\|x\|_{L,\omega} \le\alpha^{-d}\|\eta_b\|_{\omega,\operatorname{det}}^k\|x\|_{L,\omega}. \]
	Since this holds for every $\alpha\in(0,1)$,
	\[ \|\iota_A\|_\omega\le\inf_{\alpha\in(0,1)}\alpha^{-d}=1. \]

	{The tensor-slope formula \eqref{eq:tensor-slope} and $\operatorname{\widehat{\mu}}(\operatorname{det}(\overline A))=a\operatorname{\widehat{\mu}}(\overline A)$ prove \eqref{eq:arithmetic-specht-slope}.} For a $K$-basis $b=(b_1,\ldots,b_a)$ of $A$, the induced isomorphism
	\[ b^{\otimes d}:(K^a)^{\otimes d}\longrightarrow A^{\otimes d} \]
	maps $L_{a,d}^{\,S}$ onto $\operatorname{Im}(\iota_A)$ and intertwines the special linear group actions. Hence \eqref{eq:specht-invariant-space} gives $\operatorname{Im}(\iota_A)=(A^{\otimes d})^{\operatorname{SL}(A)}$. \Cref{prop:specht-absolute-irreducibility} gives absolute irreducibility of the permutation action; tensoring with the determinant line preserves this property. Permutations act isometrically at Archimedean places; their integral matrices and integral inverses preserve the non-Archimedean supremum norms. This proves all assertions.
\end{proof}

Let $W$ be a finite-dimensional $K$-vector space and let $V\subseteq W$ have dimension $s\in\mathbb{Z}_{\ge1}$. We call the line $\bigwedge\nolimits^sV\subseteq\bigwedge\nolimits^sW$ its \emph{Pl\"ucker point}. For a reductive linear algebraic $K$-group acting on $W$, we call this point \emph{geometrically semistable} if some positive-degree homogeneous invariant polynomial over an algebraic closure is nonzero on a nonzero vector of that line. We call the set of vectors whose orbit closure contains zero the \emph{null cone}. In characteristic zero, the invariant criterion identifies geometric semistability with exclusion from the null cone. These definitions are invariant under field extension; see \cite[Section 2]{bgow} and \cite[Section 5.2]{cm}. Let $z\in\bigwedge\nolimits^s W\setminus\{0\}$. Suppose that the line $Kz$ is geometrically semistable. Then a positive-degree homogeneous invariant polynomial that does not vanish at $z$ may be chosen with coefficients in $K$. Indeed, in each positive homogeneous degree the invariance identities are $K$-linear equations in the polynomial coefficients. Their solution space therefore commutes with extension of fields. If every positive-degree homogeneous invariant polynomial over $K$ vanished at $z$, then every linear combination of these polynomials over an algebraic closure would also vanish at $z$, contradicting geometric semistability.

\begin{theorem}\label{thm:arithmetic-git}
	Let $N\in\mathbb{Z}_{\ge1}$, let $(\overline E_j)_{j=1}^N$ be positive-rank Hermitian adelic bundles on $S$, and equip a nonzero subspace
	\[ i:V\hookrightarrow W:=\bigotimes_{j=1}^N E_j \]
	with induced norms. If its Pl\"ucker point is geometrically semistable under $\prod_{j=1}^N\operatorname{SL}(E_j)$, then
	\[ \operatorname{\widehat{\mu}}(\overline V) \le\sum_{j=1}^N\operatorname{\widehat{\mu}}(\overline E_j). \]
\end{theorem}
\begin{proof}
	Put $s:=\operatorname{dim}_K(V)$ and $a_j:=\operatorname{dim}_K(E_j)$ for $j\in\{1,\ldots,N\}$. For $d\in\mathbb{Z}_{\ge1}$ divisible by $s$ and every $a_j$, define
	\[
		\begin{aligned}
			B_d & :=(\operatorname{det}(\overline V))^{\otimes(d/s)}
			\otimes\overline L_{s,d}^{\,S},                          \\
			D_d & :=\bigotimes_{j=1}^N
			\left((\operatorname{det}(\overline E_j^\vee))^{\otimes(d/a_j)}
			\otimes\overline L_{a_j,d}^{\,S}\right).
		\end{aligned}
	\]
	The natural map $D_d\to(W^\vee)^{\otimes d}$ is contracting at all places by \Cref{lem:arithmetic-specht}. Define
	\begin{equation}\label{eq:arithmetic-contraction}
		\begin{gathered}
			\Psi_d:B_d\longrightarrow D_d^\vee,\\
			\forall\,(v,u)\in B_d\times D_d,\quad
			\Psi_d(v)(u):=\langle u,i^{\otimes d}(\iota_V(v))\rangle.
		\end{gathered}
	\end{equation}
	The brackets denote algebraic evaluation. Write $j_d:D_d\to(W^\vee)^{\otimes d}$ for the map just constructed. The canonical tensor-duality identification gives
	\[ \Psi_d=j_d^\vee\circ i^{\otimes d}\circ\iota_V: B_d\longrightarrow V^{\otimes d} \longrightarrow W^{\otimes d}\longrightarrow D_d^\vee. \]
	At an Archimedean place these are Hilbert--Schmidt tensor norms; at a non-Archimedean place they are the ultrametric tensor norms of \Cref{subsec:bundle-operations}. Tensoring contractions and passing to duals preserve the contraction property. The inclusion $i$ has operator norm one for the induced norms, and \Cref{lem:arithmetic-specht} applies to both $j_d$ and $\iota_V$. This factorization proves
	\begin{equation}\label{eq:arithmetic-contraction-norm}
		\forall\,\omega\in\Omega,\quad \|\Psi_d\|_\omega\le1.
	\end{equation}

	Choose $\eta\in\operatorname{det}(V)\setminus\{0\}$. By geometric semistability and the preceding argument, there is a positive-degree homogeneous invariant polynomial $p$ on $\bigwedge\nolimits^sW$, with coefficients in $K$, such that $p((\bigwedge\nolimits^s i)(\eta))\ne0$. Define
	\[ \forall\,A\in\operatorname{Hom}(V,W),\quad q(A):=p((\bigwedge\nolimits^s A)(\eta)). \]
	This polynomial is invariant under the action $A\mapsto(\bigotimes_{j=1}^N g_j)Ah^{-1}$ of $\prod_{j=1}^N\operatorname{SL}(E_j)\times\operatorname{SL}(V)$, and $q(i)\ne0$. Put $m:=\operatorname{deg}(p)$, $P:=\prod_{j=1}^N a_j$, and $d_0:=smP$. Then $d_0>0$, $s\mid d_0$, and $a_j\mid d_0$ for every $j$.

	For $\ell\in\mathbb{Z}_{\ge1}$ and $d:=\ell d_0$, the symmetric tensor representative of $q^{P\ell}$ belongs, after regrouping, to
	\[ \left(\bigotimes_{j=1}^N(E_j^\vee)^{\otimes d}\otimes V^{\otimes d} \right)^{\prod_{j=1}^N\operatorname{SL}(E_j)\times\operatorname{SL}(V)} =D_d\otimes B_d. \]
	Here we identify the auxiliary spaces with their images under the maps of \Cref{lem:arithmetic-specht}. In characteristic zero, a homogeneous polynomial has a symmetric tensor representative obtained by averaging over tensor positions. The averaging commutes with the linear group action, so the representative of $q^{P\ell}$ is invariant. The equality above then follows by taking invariants separately in each factor and using \eqref{eq:specht-invariant-space} over $K$. This representative is used only to detect nonvanishing; no local norm estimate for the averaging map is required. Write it as $\sum_{\lambda\in I}u_\lambda\otimes v_\lambda$ for a finite set $I$, where $(u_\lambda,v_\lambda)\in D_d\times B_d$ for every $\lambda\in I$. Then
	\[ q(i)^{P\ell}=\sum_{\lambda\in I} \Psi_d(v_\lambda)(u_\lambda)\ne0. \]
	Thus $\Psi_d$ is a nonzero $K$-linear map.

	The source $B_d$ and the target $D_d^\vee$ are semistable by \Cref{prop:finitegroup}. For the target, use the independent permutation groups on its factors: their external tensor product is absolutely irreducible. \Cref{prop:slope-comparison} and \eqref{eq:arithmetic-contraction-norm} therefore give
	\[ \operatorname{\widehat{\mu}}(B_d)=\operatorname{\widehat{\mu}}_{\min}(B_d) \le\operatorname{\widehat{\mu}}_{\max}(D_d^\vee)+h(\Psi_d) \le\operatorname{\widehat{\mu}}(D_d^\vee). \]
	Using \eqref{eq:degree-duality}, \eqref{eq:tensor-slope}, and \Cref{lem:arithmetic-specht}, we obtain
	\begin{equation}
		d\operatorname{\widehat{\mu}}(\overline V)+\operatorname{\widehat{\mu}}(\overline L_{s,d}^{\,S}) \le\operatorname{\widehat{\mu}}(D_d^\vee)=d\sum_{j=1}^N\operatorname{\widehat{\mu}}(\overline E_j) -\sum_{j=1}^N\operatorname{\widehat{\mu}}(\overline L_{a_j,d}^{\,S}).
	\end{equation}

	Consequently,
	\begin{equation}\label{eq:adelic-error}
		\operatorname{\widehat{\mu}}(\overline V) \le{}\sum_{j=1}^N\operatorname{\widehat{\mu}}(\overline E_j) -{\frac1d\left(\operatorname{\widehat{\mu}}(\overline L_{s,d}^{\,S}) +\sum_{j=1}^N\operatorname{\widehat{\mu}}(\overline L_{a_j,d}^{\,S})\right)}.
	\end{equation}
	The curve $S$ and all ranks are fixed. Let $\ell$ tend to infinity along $d=\ell d_0$ and use \Cref{prop:sublinear}. This proves the theorem.
\end{proof}

For a nonzero quotient $q:W\twoheadrightarrow Q$ with quotient norms, we define geometric semistability to mean geometric semistability of the Pl\"ucker point of $Q^\vee\subseteq W^\vee$. If $q$ is geometrically semistable, applying \Cref{thm:arithmetic-git} to this dual inclusion gives
\begin{equation}\label{eq:quotient-git}
	\operatorname{\widehat{\mu}}(\overline Q) \ge\sum_{j=1}^N\operatorname{\widehat{\mu}}(\overline E_j).
\end{equation}
This is the Hermitian quotient inequality in \cite[Remark 5.3.2]{cm} with its error term removed. If $\nu(\Omega_\infty)=0$, then \eqref{eq:adelic-error} already has no error term for each auxiliary degree.

\section{From the geometric inequality to arbitrary factors} \label{sec:allrank}
We use the instability-filtration argument of \cite[Theorems 5.4.3 and 5.6.1]{cm}. The geometric input is now \eqref{eq:quotient-git}. We make the comparison on graded quotients explicit by a one-sided variation of filtration weights.

\subsection{Filtrations and their graded norms}
Let $A$ be a finite-dimensional nonzero $K$-vector space. We call a family $(\mathcal F^tA)_{t\in\mathbb{R}}$ a \emph{finite real filtration} if it is decreasing, left continuous, has finitely many jumps, equals $A$ for all sufficiently negative $t$, and equals zero for all sufficiently positive $t$. Define
\[ \mathcal F^{>t}A:=\sum_{u\in(t,\infty)}\mathcal F^uA,\qquad \operatorname{gr}_{\mathcal F}^t(A) :=\mathcal F^tA/\mathcal F^{>t}A, \]
\[ \Sigma(\mathcal F):=\{t\in\mathbb{R}: \operatorname{gr}_{\mathcal F}^t(A)\ne0\},\qquad \operatorname{\mathsf e}(\mathcal F):=\frac1{\operatorname{dim}_K(A)} \sum_{t\in\Sigma(\mathcal F)}t\operatorname{dim}_K(\operatorname{gr}_{\mathcal F}^t(A)). \]
We call a basis \emph{adapted} to $\mathcal F$ if every filtration step is spanned by a subset of the basis. Its basis vectors have weights given by their last nonzero filtration steps. Any two finite real filtrations admit a common adapted basis (cf. \cite[Proposition 5.1.3]{cm}). For filtrations $\mathcal F$ and $\mathcal G$ on $A$, with weights $(x_i,y_i)\in\mathbb{R}^2$ on a common adapted basis $(b_i)_{i=1}^{\operatorname{dim}_K(A)}$, define
\[ \langle\mathcal F,\mathcal G\rangle :=\frac1{\operatorname{dim}_K(A)}\sum_{i=1}^{\operatorname{dim}_K(A)}x_i y_i, \qquad \|\mathcal F\|^2:=\langle\mathcal F,\mathcal F\rangle. \]
These quantities are independent of the common adapted basis (cf. \cite[Section 5.4]{cm}). We call $\mathcal F$ \emph{normalized} if $\operatorname{\mathsf e}(\mathcal F)=0$ and \emph{trivial} if all its weights are zero. For a quotient $q:A\twoheadrightarrow Q$, define $(\mathcal F_Q)^tQ:=q(\mathcal F^tA)$.

Let $\overline E$ and $\overline F$ be positive-rank Hermitian adelic bundles. For filtrations $\mathcal P_E$ on $E$ and $\mathcal P_F$ on $F$, define their \emph{tensor filtration} $\mathcal P$ on $W:=E\otimes F$ by
\[ \forall\,\gamma\in\mathbb{R},\quad \mathcal P^\gamma W:= \sum_{\substack{(t,u)\in\Sigma(\mathcal P_E)\times\Sigma(\mathcal P_F)\\ t+u\ge\gamma}} \mathcal P_E^tE\otimes\mathcal P_F^uF. \]
Put $E_t:=\operatorname{gr}_{\mathcal P_E}^t(E)$ and $F_u:=\operatorname{gr}_{\mathcal P_F}^u(F)$, with subquotient norms. For $\gamma\in\Sigma(\mathcal P)$, there is a canonical isomorphism
\begin{equation}\label{eq:graded-tensor-map}
	C_\gamma:= \bigoplus_{\substack{(t,u)\in\Sigma(\mathcal P_E)\times\Sigma(\mathcal P_F)\\ t+u=\gamma}} E_t\otimes F_u \longrightarrow\operatorname{gr}_{\mathcal P}^\gamma(W).
\end{equation}
The algebraic decomposition is the binary case of \cite[Proposition 5.5.5]{cm}. We use its individual component maps
\[ E_t\otimes F_u\longrightarrow \operatorname{gr}_{\mathcal P}^{t+u}(W), \]
with the tensor norms on the source and the actual subquotient norms on the target. Each component map is contracting at every place. At an Archimedean place, identify the two source grades with orthogonal complements in the filtration steps; their Hilbert--Schmidt tensor product maps isometrically into $W$, and passage to the target grade is an orthogonal quotient. At a non-Archimedean place, the tensor product of the quotient norms agrees with the corresponding quotient of the tensor product norm, and tensoring isometric inclusions preserves norms; see \cite[Proposition 1.2.36 and Remark 5.5.6]{cm}. Composing with the map to a quotient grade $Q_{t+u}:=\operatorname{gr}_{\mathcal P_Q}^{t+u}(Q)$ remains contracting. The component images span each target grade, both before and after passage to $Q$.

\begin{lemma}\label{lem:tensor-hn}
	Let $\overline E$ and $\overline F$ be positive-rank Hermitian adelic bundles. Suppose that the tensor product of every pair of their nonzero Harder--Narasimhan grades is semistable. Then the tensor product of their weighted Harder--Narasimhan filtrations is the Harder--Narasimhan filtration of $\overline E\otimes\overline F$.
\end{lemma}
\begin{proof}
	Put $\mathcal P:=\mathcal H_E\otimes\mathcal H_F$ and $G_\gamma:=\operatorname{gr}_{\mathcal P}^\gamma(E\otimes F)$. For $t+u=\gamma$, every nonzero source $E_t\otimes F_u$ is semistable of slope $\gamma$, so its HN filtration at $\gamma$ is the whole source. By HN functoriality (cf. \cite[Proposition 4.3.49]{cm}), the component contraction has image contained in $\mathcal H_{G_\gamma}^{\gamma}G_\gamma$. These images span $G_\gamma$, hence
	\[ \mathcal H_{G_\gamma}^{\gamma}G_\gamma=G_\gamma, \qquad \operatorname{\widehat{\mu}}_{\min}(G_\gamma)\ge\gamma \quad(G_\gamma\ne0). \]
	In particular $\operatorname{\widehat{deg}}(G_\gamma)\ge\gamma\operatorname{dim}_K(G_\gamma)$. The sum of the left sides is $\operatorname{\widehat{deg}}(E\otimes F)$ by degree additivity. The algebraic decomposition \eqref{eq:graded-tensor-map} and the formula for the degree of a tensor product show that the right sides have the same sum. Thus every bound is an equality. Every nonzero quotient of $G_\gamma$ has slope at least $\gamma=\operatorname{\widehat{\mu}}(G_\gamma)$, so $G_\gamma$ is semistable of that slope. HN uniqueness now proves the assertion.
\end{proof}

\subsection{A first variation on quotient filtrations}
Let $\mathcal P_E$ be a finite real filtration of $E$. We define a \emph{lifting} of the Harder--Narasimhan filtrations of its grades by constructing a filtration $\mathcal R_E$ as follows: choose an adapted basis for the Harder--Narasimhan filtration of each $E_t$, lift these vectors to $\mathcal P_E^tE$, and assign each lifted vector its Harder--Narasimhan weight in $E_t$. The union of these lifted bases is a common adapted basis for $\mathcal P_E$ and $\mathcal R_E$. This is the lifting construction of \cite[Section 5.5]{cm}.

The same lifting construction applies to the grades of a quotient filtration. HN functoriality on the component maps then controls small positive perturbations of the source weights.

\begin{lemma}\label{lem:filtered-variation}
	Let $\overline E$ and $\overline F$ be positive-rank Hermitian adelic bundles, and let $\mathcal P_E$, $\mathcal P_F$ be finite real filtrations. Choose liftings $\mathcal R_E$, $\mathcal R_F$ of the HN filtrations of their grades. Suppose that each nonzero $E_t\otimes F_u$ has the tensor product HN filtration. Let $q:E\otimes F\twoheadrightarrow Q$ be a nonzero quotient with quotient norms. On common adapted bases, replace each primary weight $x$ and lifting weight $y$ by $x+\varepsilon y$, obtaining $\mathcal P_E(\varepsilon)$ and $\mathcal P_F(\varepsilon)$. Put $\mathcal P(\varepsilon):= \mathcal P_E(\varepsilon)\otimes\mathcal P_F(\varepsilon)$ and $\mathcal P:=\mathcal P(0)$. For all sufficiently small $\varepsilon>0$,
	\begin{equation}\label{eq:filtration-variation}
		\operatorname{\mathsf e}(\mathcal P(\varepsilon)_Q) \le\operatorname{\mathsf e}(\mathcal P_Q)+\varepsilon\operatorname{\widehat{\mu}}(\overline Q).
	\end{equation}
\end{lemma}
\begin{proof}
	Write $Q_\gamma:=\operatorname{gr}_{\mathcal P_Q}^\gamma(Q)$, with the norms induced from $\overline Q$. Lift the HN filtrations of these grades to a filtration $\mathcal R_Q$ of $Q$. The degree identity for HN grades and degree additivity give
	\[ \operatorname{\mathsf e}(\mathcal R_Q) =\frac{\sum_\gamma\operatorname{\widehat{deg}}(Q_\gamma)}{\operatorname{dim}_K Q} =\operatorname{\widehat{\mu}}(\overline Q). \]
	By hypothesis, the secondary filtration on $E_t\otimes F_u$ is its HN filtration. The canonical contraction $E_t\otimes F_u\to Q_{t+u}$ preserves HN thresholds by \cite[Proposition 4.3.49]{cm}.

	Let $(x_i)$ be the tensor product of the chosen source bases, with primary and secondary weights $p_i,r_i$. Choose a basis $(y_j)$ of $Q$ adapted to both $\mathcal P_Q$ and $\mathcal R_Q$, with weights $a_j,b_j$, and write $q(x_i)=\sum_j c_{ji}y_j$. Preservation of primary weights gives $p_i\le a_j$ whenever $c_{ji}\ne0$. On the grade of weight $p_i$, preservation of HN thresholds further gives $r_i\le b_j$ when $a_j=p_i$. There are only finitely many coefficients. The strict primary gaps therefore imply, for all sufficiently small $\varepsilon>0$,
	\[ c_{ji}\ne0\quad\Longrightarrow\quad p_i+\varepsilon r_i\le a_j+\varepsilon b_j. \]

	For a weighted spanning family in an $r$-dimensional space, the sum of the weights of its quotient filtration is the maximum sum of weights over $r$ members forming a basis. To see this, choose independent vectors successively in decreasing order of weight; the resulting basis is adapted to every threshold subspace. Conversely, the number of independent vectors of weight at least any threshold is bounded by the dimension of that subspace. Apply this observation to $(q(x_i))$. Each nonzero $r\times r$ minor has a nonzero term in its determinant expansion. Summing the displayed coordinate inequalities along that term bounds its total source weight by $\sum_j(a_j+\varepsilon b_j)$. Dividing by $r=\operatorname{dim}_K Q$ yields
	\[ \operatorname{\mathsf e}(\mathcal P(\varepsilon)_Q) \le\operatorname{\mathsf e}(\mathcal P_Q) +\varepsilon\operatorname{\mathsf e}(\mathcal R_Q), \]
	which proves \eqref{eq:filtration-variation}.
\end{proof}

\subsection{Rank induction}
\begin{proof}[Proof of \Cref{thm:bost}]
	We prove simultaneously that tensor products of semistable Hermitian bundles are semistable and that arbitrary Hermitian bundles have compatible tensor HN filtrations. Use strong induction on the sum $A$ of the two positive ranks. At $A=2$ both factors are lines, and the assertions follow from the tensor-slope formula. Fix $A\ge3$ and assume both assertions for every smaller rank sum.

	First take semistable $\overline E,\overline F$ of ranks $m,n$ with $m+n=A$. It suffices to show that every nonzero quotient $q:E\otimes F\twoheadrightarrow Q$ satisfies
	\begin{equation}\label{eq:semistable-quotient}
		\operatorname{\widehat{\mu}}(\overline Q)\ge\operatorname{\widehat{\mu}}(\overline E)+\operatorname{\widehat{\mu}}(\overline F).
	\end{equation}
	For a geometrically semistable quotient this is \eqref{eq:quotient-git}. Otherwise the implication $(1)\Rightarrow(2)$ of \cite[Theorem 5.4.3]{cm}, for two factors, gives normalized filtrations $\mathcal P_E,\mathcal P_F$ attaining a negative minimum of
	\[ \Theta(\mathcal T_E,\mathcal T_F):= \frac{\operatorname{\mathsf e}((\mathcal T_E\otimes\mathcal T_F)_Q)} {(\|\mathcal T_E\|^2+\|\mathcal T_F\|^2)^{1/2}} \]
	over normalized, nontrivial pairs. Positive simultaneous scaling of the weights leaves this ratio unchanged, so choose the minimum with
	\[ \|\mathcal P_E\|^2+\|\mathcal P_F\|^2=1,\qquad c:=\operatorname{\mathsf e}((\mathcal P_E\otimes\mathcal P_F)_Q)<0. \]
	A centered filtration with only one nonzero grade is trivial. Thus at least one chosen filtration has several nonzero grades, and every pair of nonzero grades satisfies $\operatorname{dim}_K(E_t)+\operatorname{dim}_K(F_u)<A$. The induction hypothesis supplies tensor-HN compatibility for these pairs. Choose grade-HN liftings $\mathcal R_E,\mathcal R_F$. Degree additivity gives
	\[ \operatorname{\mathsf e}(\mathcal R_E)=\operatorname{\widehat{\mu}}(\overline E),\qquad \operatorname{\mathsf e}(\mathcal R_F)=\operatorname{\widehat{\mu}}(\overline F). \]
	We claim that
	\begin{equation}\label{eq:nonpositive-pairing}
		\beta:=\langle\mathcal P_E,\mathcal R_E\rangle +\langle\mathcal P_F,\mathcal R_F\rangle\le0.
	\end{equation}
	To prove the first summand is nonpositive, list the weights of $\mathcal P_E$ as $t_1>\cdots>t_h$, where $h\in\mathbb{Z}_{\ge1}$, and put $U_0:=0$ and $U_i:=\mathcal P_E^{t_i}E$ for $i\in\{1,\ldots,h\}$. The lifting construction gives
	\begin{align*}
		m\langle\mathcal P_E,\mathcal R_E\rangle
		 & =\sum_{i=1}^h t_i\operatorname{\widehat{deg}}(\overline U_i/\overline U_{i-1})       \\
		 & =t_h\operatorname{\widehat{deg}}(\overline E)
		+\sum_{i=1}^{h-1}(t_i-t_{i+1})\operatorname{\widehat{deg}}(\overline U_i)               \\
		 & \le\operatorname{\widehat{\mu}}(\overline E)
		\left(t_hm+\sum_{i=1}^{h-1}(t_i-t_{i+1})\operatorname{dim}_K(U_i)\right)                \\
		 & =m\operatorname{\widehat{\mu}}(\overline E)\operatorname{\mathsf e}(\mathcal P_E)=0.
	\end{align*}
	The inequality uses semistability of $\overline E$. Applying the same argument to $F$ proves \eqref{eq:nonpositive-pairing}.

	Put $\sigma:=\operatorname{\widehat{\mu}}(\overline E)+\operatorname{\widehat{\mu}}(\overline F)$ and let $\mathcal R_E^0,\mathcal R_F^0$ be the liftings with their mean weights subtracted. Define
	\[ D:=\|\mathcal R_E^0\|^2+\|\mathcal R_F^0\|^2\ge0, \qquad N(\varepsilon):= \operatorname{\mathsf e}(\mathcal P(\varepsilon)_Q), \]
	where $\mathcal P(\varepsilon)$ uses the uncentered lifting weights as in \Cref{lem:filtered-variation}. The centered perturbations have numerator $N(\varepsilon)-\varepsilon\sigma$ and squared norm $1+2\varepsilon\beta+D\varepsilon^2$, which is positive for all sufficiently small $\varepsilon>0$. Minimality of $c$ gives
	\[ N(\varepsilon)\ge\varepsilon\sigma+ c\sqrt{1+2\varepsilon\beta+D\varepsilon^2}. \]
	Since $\beta\le0$, $D\ge0$, and $c<0$, the elementary bound $\sqrt{1+2\varepsilon\beta+D\varepsilon^2}\le1+D\varepsilon^2$ and \Cref{lem:filtered-variation} yield
	\[ c+\varepsilon\sigma+cD\varepsilon^2 \le N(\varepsilon)\le c+\varepsilon\operatorname{\widehat{\mu}}(\overline Q). \]
	Hence $\sigma-\operatorname{\widehat{\mu}}(\overline Q)\le-cD\varepsilon$ for all sufficiently small positive $\varepsilon$. If the left side were positive, choosing a sufficiently small $\varepsilon$ would contradict this inequality. This proves \eqref{eq:semistable-quotient} and therefore semistability of $\overline E\otimes\overline F$.

	For arbitrary factors of rank sum $A$, every nonzero HN grade is semistable. Their tensor products are semistable by the semistable case just proved or by induction at a smaller rank sum. \Cref{lem:tensor-hn} gives tensor-HN compatibility, completing the simultaneous induction. The largest weight of the tensor HN filtration is $\operatorname{\widehat{\mu}}_{\max}(\overline E)+\operatorname{\widehat{\mu}}_{\max}(\overline F)$, which proves \eqref{eq:main}. Iteration gives the maximal slope identity and semistability assertion for every finite family.
\end{proof}

\begin{corollary}\label{cor:hn}
	Let $\overline E$ and $\overline F$ be positive-rank Hermitian adelic bundles on $S$, and let $\mathcal H_E$ and $\mathcal H_F$ be their weighted Harder--Narasimhan filtrations. Then
	\[ \forall\,v\in\mathbb{R},\quad \mathcal H_{E\otimes F}^v(E\otimes F) =\sum_{\substack{(t,u)\in\Sigma(\mathcal H_E)\times\Sigma(\mathcal H_F)\\ t+u\ge v}} \mathcal H_E^tE\otimes\mathcal H_F^uF. \]
\end{corollary}
\begin{proof}
	The simultaneous induction in the proof of \Cref{thm:bost} establishes this identity of tensor product filtrations. Equivalently, apply \Cref{lem:tensor-hn} to the semistable products of grades provided by that theorem.
\end{proof}

\printbibliography
\end{document}